\input amstex
\documentstyle{amsppt}
\magnification = 1100
\mathsurround = 1 pt
\define\si{\smallskip\noindent}
\define\bi{\bigskip\noindent}
\loadbold
\define\ex{\text{E}}
\define\var{\text{Var}}
\define\pr{\text{P}}
\define\la{\lambda}
\define\eps{\varepsilon}

\define\a{\alpha}
\define\part{\partial}
\define\G11{G_{1,1}(n,m)}
\topmatter
\title
Counting strongly-connected, sparsely edged directed graphs.
\endtitle
\rightheadtext{Counting digraphs}
\author
Boris Pittel
\endauthor
\affil 
Ohio State University 
\endaffil
\address
Ohio State University, Columbus, Ohio, 43210, USA 
\endaddress
\email
bgp\@math.ohio-state.edu
\endemail
\thanks
Research supported in part by NSF Grant DMS-0805996. 
\endthanks
\keywords
Enumeration, directed graphs, strongly connected, asymptotics
\endkeywords
\subjclass
05C80, 60K35
\endsubjclass
\abstract
A sharp asymptotic formula for the number of strongly connected digraphs on $n$ labelled
vertices with $m$ arcs, under a condition $m-n\to\infty$, $m=O(n)$, is obtained; 
this solves a problem posed
by Wright back in $1977$. Our formula is a counterpart of a classic  asymptotic  formula, due to
Bender, Canfield and McKay, for the total number of connected undirected graphs on $n$ vertices
with $m$ edges. A key ingredient of their proof was a recurrence equation for
the connected graphs count due to Wright. No analogue of Wright's recurrence seems to
exist for digraphs. In a previous paper with Nick Wormald we rederived the BCM
formula via counting two-connected graphs among the graphs of minimum degree
$2$, at least. In this paper, using a similar embedding for directed graphs, we
find an asymptotic formula, which includes an explicit error term, for  the fraction of strongly-connected 
digraphs with parameters $m$ and $n$ among all such digraphs with positive in/out-degrees. 
\endabstract
\endtopmatter
\document
\noindent
\bi
{\bf 0. Intoduction.\/} In a pioneering paper [24], Wright described two potent
approaches 
to calculating  $c(n,n+k)$, the number of connected labeled graphs with $n$ vertices
and $n+k$ edges. One was based on a classic exponential identity expressing a bivariate
generating function of all labeled graphs through that of the connected graphs.
This identity lent itself to a recursive procedure well tailored for computing
the exponential generating functions $W_k(x)=\sum_n x^n c(n,n+k)/n!$, and finding
the asymptotic formulas for $c(n,n+k)$ for small $k$, going beyond a formula for 
the first nontrivial $c(n,n+1)$, due to Bagaev [1]. By representing a connected
graph as a connected $2$-core with a forest of trees sprouting from the core
vertices, Wright found an alternative way of computing $W_k(x)$ for small $k$,
and established a remarkable expression of $W_k(x)$ through the exponential
generating function of the rooted trees for all $k$. Neither of these approaches
could be used to get asymptotics for $k=k(n)\to\infty$.

In [25] Wright described a quadratic recurrence
equation for the numbers $c(n,\mathbreak n+k)$, and the related recurrence for
the exponential generating function $W_k(x)$.  
The recurrence relates
to the last, $(n+k)$-th, step of an edge-insertion process that begins at an empty graph on $[n]$;
at the terminal step a connected graph with $n+k$ edges is born either if the penultimate graph
is already connected and has $n+k-1$ edges, or if it consists of two components with
parameters $(n_i,n_i+k_i)$, with $n_1+n_2=n$, $k_1+k_2+1=k$, and the $(n+k)$-th edge joins the two
components. Wright  used the recurrence to obtain a sharp asymptotic formula for $c(n,n+k)$  
with  $k=o(n^{1/3})$ as $n\to\infty$, a difficult result with many ramifications.

Bollob\'as [6], [7] discovered that a leading factor in
Wright's formula [25] was an upper bound for $c(n,n+k)$ far beyond $k=o(n^{1/3})$, and
used his bound to identify sharply---for the first time--- a transition window in the  
Erd\H os--R\'enyi random edge-insertion graph process, [9]. 
Wright's results and enumerational insight were also a key ingredient in the subsequent
studies of random graphs , \L uczak [14], \L uczak et al
[16],  Janson
et al [12], Flajolet et al [10], Daud\'e and Ravelomanana [8], Pittel
and Yeum [22], Stepanov [23], to name just a few.
 
In a monumental paper [3], Bender et al  managed to use Wright's quadratic recurrence to
derive an extension of Wright's asymptotic formula for $k=O(n\ln n)$, i. e. all the way
till the number of edges is so large that $c(n,n+k)$ is asymptotic to the number of all
graphs with parameters $n$ and $n+k$. Later with Nick Wormald we found [20], [21] a more 
combinatorial, and less 
technical
way to rederive the Bender-Canfield-McKay formula. It was based on Wright's 
decomposition device. First we counted the
connected $2$-cores among all the graphs with minimum degree $2$, at least. After this step,
enumeration of the
connected graphs was reduced to enumeration of the cores with a forest of trees attached
to the core vertices. More recently van der Hofstadt and Spencer [11] found yet another 
way to obtain the BCM formula. At the heart of their approach is a breadth-first search for
connected components in the Bernoulli-type random graph $G(n,\pr(\text{edge})=p)$.
Interestingly, this algorithm had already been used by Karp [13] to obtain
 a sharp esimate of the largest strongly-connected subgraph size of the random digraph with
$\pr(\text{arc})=c/n$ in a supercritical phase $c>1$.

Wright also wrote a companion paper [26] on determination of $g(n,\mathbreak n+k)$,
the total number of strongly-connected digraphs with $n+k$ {\it arcs\/}. He was able 
to find
a qualitative description of how $g(n,n+k)$ depends on $n$, $k$, and to find
explicit formulas for $k=0,1,2$. We gather that he did not find a recurrence for
$g(n,n+k)$ that might have been used to study the case $k=k(n)\to\infty$. (We are not
aware of any such recurrence either.)

Aside from Bender et al [4] who found an elegant way to use their formula from [3]
for enumeration of weakly-connected digraphs, the problem of asymptotic behavior of
$g(n,n+k)$ has remained open.

Our goal in this paper is to obtain a sharp asymptotic formula for $g(n,n+k)$,
with essentially best remainder term---a directed counterpart of the BCM formula---
by using a proper version 
of the embedding device. 

Consider the digraphs on $n$ labelled vertices, with
$m$ arcs. Let $g(n,m)$ stand for the total number of strongly connected digraphs, with parameters
$n$ and $m$; thus $m>n$ necessarily. Suppose $r:=m-n\to\infty$ and $m=O(n)$. We will show that
$$
\multline
g(n,m)=\frac{m!}{2\pi n\var[Y]}\,\frac{(e^{\la}-1)^{2n}}{\la^{2m}}\\
\times
\frac{\left(1-\dsize\frac{\la}{e^{\la}-1}\right)^2}{
1-\dsize\frac{\la}{e^{\la}(e^{\la}-1)}}\,\exp\left(-\frac{m}{n}-\frac{\la^2}{2}\right)\\
\times \bigl(1+O(r^{-1}\ln^2 r+r^{\eps}n^{-1/2+\gamma})\bigr),
\endmultline\tag 1
$$
for any fixed $\eps>0$, $\gamma\in (0,1/2)$. Here $\la e^{\la}/(e^{\la}-1)=m/n$, and $Y$ is
a positive Poisson variable with mean $m/n$, so that $\var[Y]=
(m/n)\bigl(\la-r/n\bigr)$. In particular, if $r=O(n^{1/2-\sigma})$, 
$\sigma>0$, then
$$
g(n,m)=\frac{m!}{6\pi e\, n}\,\left(\frac{n}{2r}\right)^{2r}\,\bigl(1+O(
r^{-1}\ln^2 r)\bigr).\tag 2
$$
Neglecting $\ln^2 r$, the error term would be of order $1/r$, which we believe
is the correct order of the remainder term.
Inevitably our proofs run parallel to the argument in [20]. There are quite a few
unexpected challenges though, as the graph component notion morphs into two, harder-to-handle, dual 
notions of a sink-set and a 
source-set, the subsets 
of vertices with no arcs going outside, and no arcs coming from outside, respectively.
\si

{\bf Note.\/}  The author is well aware of a research started earlier by Xavier 
P\'erez and Nick Wormald.  According to Wormald,
(private communication), they expect a full argument to be less technical since 
their aim is
a cruder version of (1), i.e. an explicit leading factor times $1+o(1)$. 
\si

The author plans to use the results and the insights of this paper
for a sharp analysis of a phase transition
window for the directed counterpart of the Erd\H os-R\'enyi graph process.
\bi

The rest of the paper is organized as follows. In Section 1 we state 
an asymptotic estimate for the number of digraphs with given in/out degrees (Theorem 1.1), 
analogous
to those by Bender and Canfield [2], Bollob\'as [5], McKay [17], and McKay and Wormald
[18] for graphs. In a proof sketch we introduce a random matching scheme
similar to a random pairing introduced by Bollob\'as  [5] for graphs. We use this estimate
to prove an asymptotic formula for the number of digraphs with constrained in/out degrees
(Theorem 1.2),
and to formulate auxiliary bounds we use later. In Section 2, (Theorem 2.1) we prove a bound 
for the number of those digraphs without isolated cycles. In Section 3 we use Theorem 2.1
and the bounds from Section 1 to prove (Theorem 3.1) a sharp $O(n^{-1})$ bound for the 
fraction of digraphs with a
complex ``sink-set'' (``source-set'') having less than half of all arcs, but without a simple 
sink-set (source-set). This implies that the fraction of strongly-connected digraphs
differs from the fraction of digraphs without simple sink/source-sets by at most $O(n^{-1})$.
Finally, in Section 4 we determine a sharp estimate for the latter fraction,
which turns out to be $\Theta((m-n)/n)\gg n^{-1}$ (Theorem 4.2), and this result yields
the formula (1).

\bi
{\bf 1. Enumerating the digraphs with restricted in/out-degrees.\/} 
\bi
In this Section
we provide a set of estimates, both crude and sharp, for the counts of all
digraphs with (most of) their in/out-degrees being positive. 
\proclaim{Theorem 1.1} Let $\delta_1,\dots,\delta_n\ge 0$ and $\Delta_1,\dots,\Delta_n
\ge 0$ be such that 
$$
\sum_i\delta_i=\sum_i\Delta_i=m. \tag 1.1
$$
where $m\ge n$. Introduce $g(\boldsymbol\delta,\boldsymbol\Delta)$, the total number of 
 simple
digraphs with in-degrees $\delta_i$ and out-degrees $\delta_i$. 
If $D:=\max\{\max_i\delta_i,\max_j\Delta_j\}=o(m^{1/4})$, then
$$
g(\boldsymbol\delta,\boldsymbol\Delta)=m!\left(\prod_{i\in [n]}\frac{1}{\delta_i!\Delta_i!}
\right)\,F(\boldsymbol\delta,\boldsymbol\Delta),\tag 1.2
$$
where the ``fudge factor'', always $1$ at most, is given by
$$ 
F(\boldsymbol\delta,\boldsymbol\Delta)=\exp\left(-\frac{1}{m}\sum_i\delta_i\Delta_i
-\frac{1}{2m^2}\sum_i(\delta_i)_2\sum_j(\Delta_j)_2+O(D^4/m)\right);\tag 1.3
$$
($(d)_2:=d(d-1)$).
\endproclaim
A combinatorial core of the proof is a random matching scheme, similar to a 
random {\it pairing\/} model for undirected graphs introduced by Bollob\'as [5] in his
probabilistic proof of the Bender-Canfield [2] formula for the number of graphs with
given, bounded, degrees. 
Consider two copies $[m]_1$ and $[m]_2$ of
the set $[m]$, together with the partitions $[m]_1=\uplus_i I_i$, $[m]_2=
\uplus_j O_j$, $|I_i|=\delta_i$, $|O_j|=\Delta_j$. Each of $m!$ bijections $\pi:
[m]_2\to [m]_1$ determines a directed multigraph $G(\pi)$: $i\to j$ is an
arc if $\pi(\nu)=\mu$ for some $\nu\in O_j$ and $\mu\in I_i$. 
Notice that, each {\it simple\/} digraph $G$
with in-degrees $\delta_i$ and out-degrees $\Delta_j$, corresponds to
exactly $\prod_i(\delta_i!\Delta_i!)$ bijections $\pi$. We call
those bijections digraph-induced. Therefore
$F(\boldsymbol\delta,\boldsymbol\Delta)$ {\it defined\/} by (2) is the probability
that the bijection $\pi$ chosen uniformly at random from among all $m!$
bijections is digraph-induced, i.e. there is no $i$ such that $\pi(\nu)=\mu$ for
$\nu\in O_i$, $\mu\in I_i$, and there is no $(i,j)$ such that $\pi(\nu_t)=\mu_t$,
$t=1,2$, ($\nu_1\neq \nu_2$, $\mu_1\neq \mu_2$), and $\nu_t\in O_j$, $\mu_t\in O_i$.
\si

Since $F(\boldsymbol\delta,\boldsymbol\Delta)\le 1$, we see that
$$
g(\boldsymbol\delta,\boldsymbol\Delta)\le m!\left(\prod_{i\in [n]}\frac{1}{\delta_i!\Delta_i!}
\right)\tag 1.4
$$
always.

So the task is 
to show that this probability $F(\boldsymbol\delta,\boldsymbol\Delta)\le 1$ is 
given by (1.3).  The sums in the exponent 
in (1.3) are asymptotic estimates for
the expected number of loops and number of double arcs in $G(\pi)$. Our claim is a 
digraph analogue of McKay's far-reaching extension [17] (McKay and
Wormald [18]), based on Bollob\'as' pairing model, of the Bender-Canfield
formula  for the number of undirected graphs with a given degree 
sequence $d_1,\dots,d_n$, and $max_i\,d_i=o(M^{1/4})$. We omit the proof as 
it is very similar to those in [17] and [18]. 
\bi

{\bf Note 1.1.\/} Notice that the RHS of (1.4) is not necessarily an integer, but
$$
h(\boldsymbol\delta,\boldsymbol\Delta):= 
(m!)^2\left(\prod_{i\in [n]}\frac{1}{\delta_i!\Delta_i!}
\right)=\binom{m}{\boldsymbol\delta}\,\binom{m}{\boldsymbol\Delta}
$$
most certainly is! In fact, $h(\boldsymbol\delta,\boldsymbol\Delta)$ counts the
total number of $m$-long ordered sequences of insertions of single arcs into an initially
empty digraph that result in a {\it multi\/}-digraph with in-degrees 
$\boldsymbol\delta$ and out-degrees $\boldsymbol\Delta$. 

This formula follows from a bijection between that the set of those
sequences and the set of pairs of $m$-long words in the alphabet of
$n$ letters $1,\dots, n$, such that the first (second, resp.) word has
$\delta_i$ ($\Delta_i$, resp.) letters $i$. Observe also that, given
$\boldsymbol\mu=\{\mu_{j,i}\}_{j,i\in [n]}$ with $\|\boldsymbol\mu\|=m$,
the number $h(\boldsymbol\mu)$ of the sequences, that result in the multigraph
with exactly $\mu_{j,i}$ arcs $j\to i$, is obviously the multinomial coefficient
$$
h(\boldsymbol\mu)=\binom{m}{\boldsymbol\mu}.\tag 1.5
$$

On the other hand, 
$m!\,g(\boldsymbol\delta,\boldsymbol\Delta)$ is the total number of the $m$-long
sequences that result in a simple digraph with in-degrees 
$\boldsymbol\delta$ and out-degrees $\boldsymbol\Delta$. We call these special sequences
graphic. Thus $F(\boldsymbol\delta,\boldsymbol\Delta)$ is the fraction of graphic
sequences among all $m$-long insertion sequences leading to a multi-digraph with
the prescribed in/out-degrees.
\bi

McKay's formula was used by Pittel and Wormald [19] to derive a sharp estimate of $C_k(n,m)$,
the total
number of graphs on $[n]$ with $m=O(n\ln n)$ edges, and with minimum vertex degree $k$, 
at least, ($k\ge 1$). We use Theorem 1.1 to derive an asymptotic formula
for $C_{1,1}(n,m)$, the total number of digraphs on $[n]$ with $\min_i\delta_i\ge 1$,
$\min\Delta_i\ge 1$. In fact we will need a slightly more general set-up. 
\si

Given $n_1\le n$, let $\bold n=(n_1,n)$, and let
$C_{1,1}(\bold n,m)$ denote the total number of digraphs on $[n]$, with $m$ arcs, such that
$$
\delta_i\ge 1,\, (i\in [n]);\quad \Delta_j\ge 1,\,(j\in [n_1]);\quad 
\Delta_j\equiv 0,\,(j>n_1).\tag 1.6
$$
By (1.4),
$$
\aligned
C_{1,1}(\bold n,m)\le&\, m!\sum_{\boldsymbol\delta,\boldsymbol\Delta
\text{ meet }(1.6)\atop |\boldsymbol\delta|=|\boldsymbol\Delta|=m}\,
\prod_{i\in [n]}\frac{1}{\delta_i!\Delta_i!}\\
=&m!\,[x^m x_1^m]\sum_{\boldsymbol\delta,\boldsymbol\Delta\text{ meet }(1.6)}\,
\prod_{i\in [n]}\frac{x^{\delta_i}x_1^{\Delta_i}}{\delta_i!\Delta_i!}\\
=&m!\,[x^m] f(x)^n\, [x_1^m]f(x_1)^{n_1},\quad f(y):=e^y-1.
\endaligned
$$
So, for all $x>0,x_1>0$,
$$
C_{1,1}(\bold n,m)\le m!\,\frac{f(x)^n}{x^m}\,\frac{f(x_1)^{n_1}}{x_1^m},\tag 1.7
$$
and the best values of $x$ and $x_1$ are the minimum points of the first fraction and the
second fraction, i.e. $\la$ and $\la_1$, the roots of
$$
\frac{\la e^{\la}}{e^{\la}-1}=\frac{m}{n},\qquad\frac{\la_1 e^{\la_1}}{e^{\la_1}-1}=
\frac{m}{n_1}.\tag 1.8
$$
Using the Cauchy integral formula
$$
[y^a]\,f(y)^b=\frac{1}{2\pi i}\oint\limits_{|z|=\rho}\frac{f(z)^b}{z^{a+1}}\,dz,
$$
and an inequality
$$
|e^z-1|\le (e^{|z|}-1)\,\exp\left(-\frac{|z|-\text{Re }z}{2}\right),
$$
one can show easily that, in fact,
$$
[y^a]\,f(y)^b\le_b\,\frac{1}{(yb)^{1/2}}\,\frac{f(y)^b}{y^a},
\quad\forall\, y>0.
$$
(Here and elsewhere $A\le_b B$ means that $A=O(B)$ uniformly over all parameters 
that determine the
values of $A$ and $B$.) So (1.7) can be strengthened to
$$
C_{1,1}(\bold n,m)\le_b m!\,\frac{f(x)^n}{(1+nx)^{1/2}x^m}\,\cdot
\frac{f(x_1)^{n_1}}{(1+n_1x_1)^{1/2}x_1^m},\tag 1.9
$$
for all $x,x_1>0$. In particular,
$$
C_{1,1}(\bold n,m)\le_b m!\,\frac{f(x)^{2n}}{(1+nx)x^{2m}},\quad\forall\,x>0.\tag 1.10
$$
\si
The next theorem gives an asymptotically sharp formula for $C_{1,1}(\bold n,m)$.
\proclaim{Theorem 1.2\/} Let $n_1=\Theta(n)$, $r:=m-n\to\infty$, $m=O(n)$. 
Introduce $Y$ and $Y_1$, two positive Poissons, with parameters $\la$ and $\la_1$, 
i.e.
$$
\pr(Y=j)=\frac{\la^j/j!}{f(\la)},\qquad \pr(Y=j)=\frac{\la_1^j/j!}{f(\la_1)};\qquad 
(j\ge 1).
$$
Then
$$
\aligned
C_{1,1}(\bold n,m)=&\,\bigl(1+O(r^{-1}+n^{-1/2}r^{\eps})\bigr)\,m!\,\frac{f(\la)^nf(\la_1)^{n_1}}
{(\la\la_1)^m}\\
&\times\frac{e^{-\eta}}{2\pi(n\var\,[Y])^{1/2}(n_1\var\,[Y_1])^{1/2}},
\endaligned\tag 1.11
$$
for every $\eps>0$, where
$$
\eta=\frac{m}{n}+\frac{\la\la_1}{2}.\tag 1.12
$$
\endproclaim
\bi

The idea of the proof is that, introducing the independent copies $Y^1,\dots,Y^n$ of $Y$
and $Y_1^1,\dots,Y_1^{n_1}$ and denoting $\bold Y=(Y^1,\dots,Y^n)$, $\bold Y_1
=(Y_1^1,\dots,Y_1^{n_1}, n-n_1\text{ zeroes})$, we can rewrite (1.2) in the following, suggestive, way:
$$
\align
C_{1,1}(\bold n,m)=&\,m!\,\frac{f(\la)^nf(\la_1)^{n_1}}{\la^m\la_1^{m_1}}\\
&\cdot [x^mx_1^{m_1}]\sum_{\boldsymbol\delta,\boldsymbol\Delta\text{ meet }(1.6)}
F(\boldsymbol\delta,\boldsymbol\Delta)\prod_{i\in [n]}\frac{(\la x)^{\delta_i}/\delta_i!}{f(\la)}
\cdot\frac{(\la_1x_1)^{\Delta_i}/\Delta_i!}{f(\la_1)}\\
=&\,m!\,\frac{f(\la)^nf(\la_1)^{n_1}}{\la^m\la_1^{m_1}}\,\ex\left[F(\bold Y,\bold Y_1)
\bold 1_{\{\|\bold Y\|=\|\bold Y_1\|=m\}}\right]\\
=&\,m!\,\frac{f(\la)^nf(\la_1)^{n_1}}{\la^m\la_1^{m_1}}\,\ex\biggl[F(\bold Y,\bold Y_1)\,
\boldkey|\,
\|\bold Y\|=\|\bold Y_1\|=m\biggr]\tag 1.13\\
&\times\pr(\|\bold Y\|=m)\,\,\pr(\|\bold Y_1\|=m).\tag 1.14
\endalign
$$
Here
$$
\ex\left[F(\bold Y,\bold Y_1)
\bold 1_{\{\|\bold Y\|=\|\bold Y_1\|=m\}}\right]=\ex\bigl[F(\bold{Z},\bold{Z}_1)\bigr],\tag 1.15
$$
where, denoting $\bold i=(i_1,\dots,i_n)$, $\bold j=(j_1,\dots,j_n)$,
$$
\pr(\bold Z=\bold i)=\frac{\prod_{s=1}^n1/i_s!}{\sum\limits_{\bold j:\,\|\bold j\|=m}\prod_{t=1}^n1/j_t!},
\quad (\|\bold i\|=m),\tag 1.16
$$
and, denoting $\bold i=(i_1,\dots,i_{n_1})$, $\bold j=(j_1,\dots,j_{n_1})$,
$$
\pr(\bold Z_1=\bold i)=\frac{\prod_{s=1}^{n_1}1/i_s!}{\sum\limits_{\bold j:\,\|\bold j\|=m}\prod_{t=1}^{n_1}1/j_t!},
\quad (\|\bold i\|=m).\tag 1.17
$$

The next step is to prove two similar cases of a local limit theorem
$$
\aligned
\pr(\|\bold Y\|=m)=&\,\frac{1+O(r^{-1})}{(2\pi n\var[Y])^{1/2}},\\ 
\pr(\|\bold Y_1\|=m)=&\,\frac{1+O(r^{-1})}{(2\pi n_1\var[Y_1])^{1/2}},
\endaligned\tag 1.18
$$
with the remainder term $O(r^{-1})$.

The last step
is to prove that on the
event $\{|\bold Y|=|\bold Y_1|=m\}$, in probability  $F(\bold Y,\bold Y_1)$ is within 
$1+o(1)$ factor from the RHS of (1.3), where $\delta_i$, $(\delta_i)_2$, ($i\in
[n]$),
are replaced with $\ex[Y]$ and $\ex[(Y)_2]$, and $\Delta_i$, $(\Delta_i)_2$, ($i\in [n_1]$),
are
replaced with $\ex[Y_1]$ and $\ex[(Y_1)_2]$. We omit the technical details
as the argument is a natural modification of 
the proof of a corresponding formula in [19]
for the total number of undirected graphs of mindegree $k\ge 1$, at least. We mention
only that in (1.11) $r^{\eps}$ comes from an observation that, with sufficiently
high probability, $\max\{\max_iY^i,\,\max_jY_1^j\}<r^{\eps}$, and that $r^{-1}$ comes
from the remainder term in (1.18).\qed
\bi

{\bf Note 1.2.\/} Implicit in the theorem 1.2 and its proof is the following:  for $r\to\infty$,
$m=O(n)$,
$$
[x^m]\,(e^x-1)^n=\frac{(e^{\la}-1)^n}{\la^{2m}}\cdot\frac{1+O(r^{-1})}{(2\pi n\var[Y])^{1/2}}.\tag 1.19
$$
\bi
\proclaim{Corollary 1.3} Let $r\to\infty$, $m=O(n)$. Then
$$
\aligned
C_{1,1}(n,m)=&\,\bigl(1+O(r^{-1}+n^{-1/2}r^{\eps})\bigr)\,m!\,\frac{f(\la)^{2n}}
{\la^{2m}}
\cdot\frac{e^{-\eta}}{2\pi n\var\,[Y]},\\
\eta=\eta(n,m):=&\,\frac{m}{n}+\frac{\la^2}{2},
\endaligned\tag 1.20
$$
for every $\eps>0$.
\endproclaim
\bi

Since $\var[Y]=\Theta(\la)$, $\var[Y_1]=\Theta(\la_1)$, the bound (1.10) for
$C_{1,1}(\bold n,m)$ with
$x=\la$, $x_1=\la_1$ contains all the key factors in the sharp estimate (1.15), except 
for $e^{-\eta}$. 
\bi

{\bf Note 1.3\/} Recalling the note 1.1, $m!\, C_{1,1}(n,m)$ is the total number
of $m$-long insertion sequences $\bold s$ resulting in a digraph with 
non-zero in/out-degrees. Let $h_{1,1}(n,m)$ denote the total number of those
sequences resulting in multi-digraphs with non-zero in/out-degrees. 
Dropping the fudge factor $F$, we have a counterpart of (1.13)-(1.15):
$$
h_{1,1}(n,m)=\,(m!)^2\,\frac{f(\la)^{2n}}{\la^{2m}}\,
\times\pr(\|\bold Y\|=m)\,\,\pr(\|\bold Y_1\|=m).\tag 1.21
$$
Consequently,
$$
\frac{m!\,C_{1,1}(n,m)}{h_{1,1}(n,m)}=\ex\bigl[F(\bold{Z},\bold{Z}_1)\bigr],\tag 1.22
$$
where $\bold Z,\,\bold Z_1$ are defined in (1.16)-(1.17), with $n_1=n$ this time.
Now $\ex\bigl[F(\bold{Z},\bold{Z}_1)\bigr]$ is asymptotically equivalent to
$e^{-\eta}$, which is positive if $r=O(n)$. Thus, for $r=O(n)$  the graphic $m$-long 
sequences constitute an asymptotically {\it positive\/} fraction among all $m$-long
sequences of arc insertions. Moreover, from the proof of the theorem 1.2 it
follows that, for $r=O(n)$,
$$
\aligned
\bigl|\bigl\{\bold s\,:\,\bigl|F(\boldsymbol\delta(\bold s),
\boldsymbol\Delta(\bold s))-e^{-\eta}\bigr|\le& r^{\eps}n^{-1/2+\gamma}\bigr\}\bigr|\\
\ge& h_{1,1}(n,m)\bigl(1-O(e^{-r^{c\eps}n^{2\gamma}})\bigr),
\endaligned\tag 1.22
$$ 
for all $\eps>0$, $\gamma\in (0,1/2)$, and some absolute constant $c>0$.
\bi
{\bf Section 2. Counting digraphs with restricted in/out-degrees that have no
isolated cycles.\/}
\bi
Recall that, given $n_1\le n$, $C_{1,1}(\bold n,m)$ denotes 
the total number of digraphs with in-degrees of
all vertices and 
outdegrees of the vertices in $[n_1]$ each being $1$, at least. Introduce
$C_{1,1}^0(\bold n,m)$, the total number of these digraphs without isolated
cycles. How much smaller is $C_{1,1}^0(\bold n,m)$ compared with $C_{1,1}(\bold n,m)$?
\proclaim{Theorem 2.1} If $(\bold n=(n_1,n_2),m)$ meet the conditions of Theorem 1.2,
and $r_1:=m-n_1=O(r)$, then
$$
C_{1,1}^0(\bold n,m)\le_b\frac{r}{n}\,C_{1,1}(\bold n,m),\quad r=m-n.\tag 2.1
$$
\endproclaim

\noindent
{\bf Proof of Theorem 2.1.\/} Obviously, we may and will assume that $r=o(n)$. If so,
in (1.10) $\eta=O(1)$, see (1.11).
Suppose a digraph is chosen uniformly at random
from among all $C_{1,1}(\bold n,m)$ digraphs. Let $X$ denote the total number of
isolated cycles in the random digraph. We need to show that $\pr(X=0)=O(r/n)$.
\si

Since the exponential generating function of directed cycles is 
$$
\sum_{\ell\ge 2}\frac{x^{\ell}}{\ell}=\ln (1-x)^{-1}-x,
$$
the binomial moments of $X$ are given by 
$$
\align 
\ex\left[\binom{X}{k}\right]=&\,\frac{1}{k!}\sum_{a< n_1}\binom{n_1}{a}
\frac{C_{1,1}(\bold n(a),m-a)}{C_{1,1}(\bold n,m)}\cdot a!\,[x^a]
\left(\sum_{\ell\ge 2}\frac{x^{\ell}}{\ell}\right)^k\\
=&\,\frac{1}{k!}\sum_{a< n_1}R_{\bold n,m}(a)\,[x^a]
\left(\sum_{\ell\ge 2}\frac{x^{\ell}}{\ell}\right)^k,\\
R_{\bold n,m}(a):=&\,\frac{n_1!C_{1,1}(\bold n(a),m-a)}{(n_1-a)!C_{1,1}(\bold n,m)},\quad 
(\bold n(a)=(n_1-a,n_2).
\endalign
$$
{\it Explanation\/}: The binomial moment is the expected number of $k$-long unordered
tuples of isolated cycles. To generate those tuples in all digraphs
in question we choose a vertex set of a generic cardinality $a$ from $[n_1]$ in
$\binom{n_1}{a}$ ways, then form an unordered family of $k$ disjoint cycles on these
$a$ vertices, in $(a!/k!)[x^a]\bigl(\sum_{\ell}x^{\ell}/\ell\bigr)^k$ ways, and finally
select an admissible complementary digraph on the remaining $n-a$ vertices with $m-a$ arcs,
in $C_{1,1}(\bold n(a),m-a)$ ways. Furthermore, an admissible $a$ also satisfies $m-a\le 
(n-a)^2$, which implies that $n-a\ge \sqrt r$.
\si

Consequently, by using an inversion formula
$$
\pr(X=j)=\sum_{k\ge j}(-1)^{k-j}\binom{k}{j}\ex\left[\binom{X}{k}\right],
$$
we obtain
$$
\pr(X=0)=\sum_{k\ge 0}(-1)^k\ex\left[\binom{X}{k}\right]=\sum_{a}R_{\bold n,m}(a)\,[x^a]h(x),
\tag 2.2
$$
where
$$
h(x)=\exp\left(-\sum_{\ell\ge 2}\frac{x^{\ell}}{\ell}\right)=(1-x)e^x.
$$
So
$$
[x^a]h(x)=\frac{1}{a!}-\frac{1}{(a-1)!},\quad a\ge 1,\quad [x^0]h(x)=1.
$$
Let us find a sharp asymptotic formula for the sum in (2.2). Using (1.10) for $C_{1,1}
(\bold n,m)$ and (1.8), with $x=\la$, $x_1=\la_1$, for $C_{1,1}(\bold n(a),m-a)$, we
bound
$$
\align
R_{\bold n,m}(a)\le_b&\,\frac{(n_1)_a}{(m)_a}\cdot\left(\frac{\la\la_1}
{f(\la)f(\la_1)}\right)^a\cdot\sqrt{\frac{nn_1}{(n-a)(n_1-a)}}\\
\le&\,\frac{n}{n_1-a}\,\sigma^a.
\endalign
$$
Here, by (1.7),
$$
\sigma=\frac{n_1}{m}\cdot\frac{\la\la_1}{f(\la)f(\la_1)}=\frac{\la e^{-\la_1}}{e^{\la}-1}
<1.\tag 2.3
$$
Introduce 
$$
A_n=3\left\lceil \frac{\ln(n/r)}{\ln\ln (n/r)}\right\rceil. 
$$
Consider $a\ge A_n$:
$$
\aligned
\left|\sum_{a\ge A_n}R_{\bold n,m}(a)\,[x^a]h(x)\right|\le_b&\, \sum_{A_n\le a\le n_1/2}
\frac{1}{a!}\,+\,\frac{n}{\lfloor n_1/2\rfloor !}\\
\le_b&\frac{1}{A_n!}+\frac{n}{\lfloor n_1/2\rfloor !}\\
\le_b&\,\left(\frac{r}{n}\right)^2.
\endaligned\tag 2.4
$$

For $1\le a\le A_n$ we need a sharp estimate  $R_{\bold n,m}(a)$, within a factor $1+O(r/n)$.
Observe upfront that the asymptotic estimate (1.10) used separately for \linebreak $C_{1,1}
(\bold n,m)$ and $C_{1,1}(\bold n(a),m-a)$ would not work, as the remainder term 
$O(r^{-1}+n^{-1/2+\eps})$ is too big. (We already confronted a similar issue in [20];
the hurdle is higher this time, as in [20] we were content to have an error
term of order $o\bigl(r/n)^{1/2}\bigr)$, rather than  $O(r/n)$. For the next theorem,
at a similar point we will even need $o(r/n)$.)
\si

Instead of (1.10), we use the exact formulas (1.12)-(1.13) for both 
$C_{1,1}(\bold n,m)$ and $C_{1,1}(\bold n(a),\mathbreak m-a)$, using $\lambda$ and
$\lambda_1$, the roots of (1.7), for both $C_{1,1}(\bold n,m)$ and
$C_{1,1}(\bold n(a),m-a)$. So
$$
\multline
R_{\bold n,m}(a)=\frac{(n_1)_a}{(m)_a}\left(\frac{\la\la_1}{f(\la)f(\la_1)}\right)^a
\cdot\frac{E_{\bold n(a),m-a}}{E_{\bold n,m}}\\
\times\frac{\pr\bigl(\sum_{i=1}^{n-a}Y^i=m-a\bigr)}{\pr\bigl(\sum_{i=1}^nY^i=m\bigr)}
\cdot\frac{\pr\bigl(\sum_{j=1}^{n_1-a}Y^j_1=m-a\bigr)}{\pr\bigl(\sum_{j=1}^{n_1}Y_1^j=m\bigr)},
\endmultline\tag 2.5
$$
where $E_{\bold n,m}=\left.E_{\bold n(a),m-a}\right|_{a=0}$, 
(see (1.14)-(1.16))
$$
\multline
E_{\bold n(a),m-a}=\ex\biggl[\bold 1_{B(a)}\exp\biggl(-\frac{1}{m-a}\sum_{i=1}^{n_1-a}Z^i(a)Z_1^i(a)\\
-
\frac{1}{2(m-a)^2}\sum_{i\notin (n_1-a,n_1]}(Z^i(a))_2\sum_{j=1}^{n_1-a}(Z_1^j(a))_2\\
+O\bigl(m^{-1}\max_{i,j}(Z^i(a)+Z_1^j(a))^4)\bigr)\biggr)\biggr]+O\bigl(\pr(B(a)^c)
\bigr),
\endmultline\tag 2.6
$$
and 
$$
B(a):=\{\max_{i,j}(Z^i(a)+Z_1^j(a))\le \omega\},
$$
with $\omega=o(m^{-1/4})$ to be specified shortly.
Our notation
emphasizes dependence of $Z$'s on 
$a$: for instance, $Z^1(a),\dots,Z^{n-a}(a)$ are occupancy numbers in the
uniformly random allocation of $m-a$ distinguishable balls among $n-a$ boxes, subject
to the condition ``no box is empty''.
\si

First, let us dispense with the ratios of local probabilities. In [20] the 
following estimate was proved. Let $Y^1,\dots,Y^n$ be 
independent Poissons distributed on
$\{2,3,\dots\}$, such that $\ex[Y^j]=2m/n$. Then, for $a=o(n)$,
$$
\frac{\pr\bigl(\sum_{i=1}^{n-a}Y^i=2(m-a)\bigr)}{\pr\bigl(\sum_{i=1}^nY^i=2m\bigr)}=
1+O(an^{-1}+a^2rn^{-2}),\quad r:=2(m-n).
$$
No real changes are needed to show that for our Poissons $Y^1,\dots,Y^n$ 
(distributed on $\{1,2,\dots\}$),
$$
\frac{\pr\bigl(\sum_{i=1}^{n-a}Y^i=m-a\bigr)}{\pr\bigl(\sum_{i=1}^nY^i=m\bigr)}=
1+O(an^{-1}+a^2rn^{-2}),\quad r:=m-n,\tag 2.7
$$
and likewise
$$
\frac{\pr\bigl(\sum_{j=1}^{n_1-a}Y^j_1=m-a\bigr)}{\pr\bigl(\sum_{j=1}^{n_1}Y_1^j=m\bigr)}=
1+O(an^{-1}+a^2r_1n^{-2}),\quad r_1=m-n_1.
$$
Since $r_1=O(r)$ we conclude that
$$
\aligned
\frac{\pr\bigl(\sum_{i=1}^{n-a}Y^i=m-a\bigr)}{\pr\bigl(\sum_{i=1}^nY^i=m\bigr)}
&\cdot\frac{\pr\bigl(\sum_{j=1}^{n_1-a}Y^j_1=m-a\bigr)}
{\pr\bigl(\sum_{j=1}^{n_1}Y_1^j=m\bigr)}\\
=&\,1+O(an^{-1}+a^2rn^{-2}).
\endaligned\tag 2.8
$$
Next we focus on $E_{\bold n(a),m-a}/E_{\bold n,m}$. Let us look at $\max_{i\notin (n_1-a,n_1]}
Z^i(a)$. For  $j\ge 1$, by (1.15),
$$
\align
\pr\bigl(Z(a)=j\bigr)=&\,\frac{1}{j!}\cdot\frac{\sum\limits_{\bold i>\bold 0:\|\bold i\|=m-a-j}\prod_{s=2}^{n-a}1/i_s!}
{\sum\limits_{\bold j>\bold 0:\,\|\bold j\|=m-a}\prod_{t=1}^{n-a}1/j_t!},\\
=&\,\frac{1}{j!}\,\frac{[x^{m-a-j}](e^x-1)^{n-a-1}}{[x^{m-a}](e^x-1)^{n-a}}\\
\le_b&\,\frac{1}{j!}\,\frac{(e^{\la}-1)^{n-a-1}/\bigl[\la^{m-a-j}\sqrt{(n-a-1)\la}\,\bigr]}
{(e^{\la}-1)^{n-a}/\bigl[\la^{m-a}\sqrt{(n-a)\la}\,\bigr]}\\
\le_b&\,\frac{1}{j!}\cdot\frac{\la^j}{e^{\la}-1}.
\endalign
$$
As $\la\le 2r/n$, we see then that $\pr(Z(a)>j)\le_b\, (2r/n)^j$. Therefore, picking $\eps>0$,
$$
\align
\pr(\max_iZ^i(a)>r^{\eps})\le& n\left(\frac{2r}{n}\right)^{r^{\eps}}\\
\le_b&\frac{r}{n}\exp\bigl[\ln n-r^{\eps}\ln(n/(2r))\bigr]\\
\le&\,\frac{r}{n}\left(\frac{r}{n}\right)^{r^{\eps}/2},\tag 2.9
\endalign
$$
since, for $r\to\infty$ and $r=o(n)$,
$$
r^{\eps}\ln(n/r)\gg \ln n.
$$
And on the event $\{\max_iZ^i(a)\le r^{\eps}\}$, we have 
$$
\frac{\max_i Z^i(a)^4}{m}\le_b \frac{r^{4\eps}}{n}= \frac{r}{n}\,r^{4\eps-1}\ll \frac{r}{n},
\tag 2.10
$$
provided that $\eps<1/4$. Obviously the counterparts of (2.9), (2.10)  hold for
$\max_jZ^j_1(a)$ as well.

So, setting $\omega=r^{\eps}$ in (2.6),
$$
\multline
E_{\bold n(a),m-a}=O\bigl((r/n)^{r^{\eps}/2+1}\bigr)+\bigl(1+O(r^{4\eps}/n)\bigr)\\
\times\ex\biggl[
\exp\biggl(-\frac{1}{m-a}\sum_{i=1}^{n_1-a}Z^i(a)Z_1^i(a)\\
-
\frac{1}{2(m-a)^2}\sum_{i\notin (n_1-a,n_1]}(Z^i(a))_2\sum_{j=1}^{n_1-a}(Z_1^j(a))_2\biggr)
\biggr],
\endmultline
$$
and, by $\ex[e^U]\ge e^{\ex[U]}$, the last expected value is bounded below by
$$
\exp\biggl(-\frac{n_1-a}{m-a}\,\ex[Z(a)]\,\ex[Z_1(a)]-\frac{(n-a)(n_1-a)}
{2(m-a)^2}\,\ex[(Z(a))_2]\,\ex[(Z_1(a))_2]\biggr),
$$
which is bounded away from $0$. Therefore
$$
\align
E_{\bold n(a),m-a}=&\,\bigl(1+O((r/n)^{r^{\eps}/2+1}+r^{4\eps}/n)\bigr)E_{\bold n(a),m-a}^*,\\
E_{\bold n(a),m-a}^*:=&\,\ex\biggl[
\exp\biggl(-\frac{1}{m-a}\sum_{i=1}^{n_1-a}Z^i(a)Z_1^i(a)\\
&-
\frac{1}{2(m-a)^2}\sum_{i\notin (n_1-a,n_1]}(Z^i(a))_2\sum_{j=1}^{n_1-a}(Z_1^j(a))_2\biggr)
\biggr],
\endalign
$$
It remains to consider $E_{\bold n(a),m-a}^*/E_{\bold n(0),m}^*$.
Notice that, conditioned on the event 
$$
A=\{Z^i(0)\equiv 1,\,i\in (n_1-a,n_1]\}\cap\{Z^i_1(0)\equiv 1,\,i\in (n_1-a,n_1]\},
$$
$\{Z^i(0),\,Z^j_1(0)\,;\,i\notin (n_1-a,n_1],\,j\le n_1-a\}$ has the same distribution as 
$\{Z^i(a),
\mathbreak\,Z^j_1(a)\,:\,i\notin (n_1-a,n_1] ,\,j\le n_1-a\}$,
and 
$$
\pr(A)\ge\, 1-a\bigl(\pr(Z(0)>1)+\pr(Z_1(0)>0)\bigr)=1-O(a\la)=1-O(ar/n).
$$
Therefore
$$
\align
E_{\bold n(a),m-a}^*=&\,\bigl(1+O(ar/n)\bigr)E_{\bold n(a),m-a}^{**},\\
E_{\bold n(a),m-a}^{**}:=&\,\ex\biggl[
\exp\biggl(-\frac{1}{m-a}\sum_{i=1}^{n_1-a}Z^i(0)Z_1^i(0)\\
&-
\frac{1}{2(m-a)^2}\sum_{i\notin (n_1-a,n_1]}(Z^i(0))_2\sum_{j=1}^{n_1-a}(Z_1^j(0))_2\biggr)
\biggr].\tag 2.11
\endalign
$$
The contribution to the expectation $E_{\bold n(a),m-a}^{**}$ from the random outcomes with
$\max_{i,j}(Z^i(0)+Z_1^j(0))>r^{\eps}$ is of order 
$$
\pr(\max_{i,j}(Z^i(0)+Z_1^j(0))>r^{\eps})\le (r/n)^{r^{\eps}/2+1}.
$$
If $\max_{i,j}(Z^i(0)+Z_1^j(0))\le r^{\eps}$, then the difference between the random exponents
in (2.11) for $a=0$ and $a>0$ is, by simple algebra, of order $ar^{4\eps}/n$. So
 $E_{\bold n(a),m-a}^{**}$ is within the multiplicative factor
$$
\bigl(1+O((r/n)^{r^{\eps}/2+1}\bigr)\,\left[\exp\bigl(O(r^{4\eps}/n)\bigr)\right]^a
$$
away from $E_{\bold n(0),m}^*$. Collecting the pieces
we obtain 
$$
\frac{E_{\bold n(a),m-a}}{E_{\bold n,m}}=\bigl(1+O(ar/n)\bigr)\,
\left[\exp\bigl(O(r^{4\eps}/n)\bigr)\right]^a,\tag 2.12
$$
uniformly for $1\le a\le A_n$.
\si

Combining (2.5), (2.8) and (2.12), we conclude: uniformly for $1\le a\le A_n$,
$$
\multline
R_{\bold n,m}(a)=\,\frac{(n_1)_a}{(m)_a}\,\left(\frac{\la\la_1}{f(\la)f(\la_1)}\,
e^{O(r^{4\eps}/n)}\right)^a
\,(1+O(ar/n)),\\
=\,\left(\frac{\la\la_1}{(m/n_1)f(\la)f(\la_1)}\,e^{O(r^{4\eps}/n)}
\right)^a\cdot\exp\left(-\frac{a^2}{2n_1}+\frac{a^2}{2m}\right)
\,(1+O(ar/n))\\
=\,\bigl[\sigma e^{O(r^{4\eps}/n)}\bigr]^a (1+O(ar/n)),
\endmultline
$$
see (2.3) for the definition of $\sigma$. Since $\sigma=1-\Theta(r/n)$, and $\eps<1/4$,
$$
\sigma e^{O(r^{4\eps}/n)}=1-\Theta(r/n).
$$
Therefore, invoking the inversion formula (2.2) and 
(2.4), {\it and\/} using 
$$
\frac{1}{a!}-\frac{1}{(a-1)!}\le 0,\quad\forall\,a\ge 1,
$$
we see that $\pr(X=0)$ is bounded above and below by
$$
\align
&\,1+\sum_{a=1}^{A_n}\left(\frac{1}{a!}-\frac{1}{(a-1)!}\right)
\bigl[\sigma e^{O(r^{4\eps}/n)}\bigr]^a(1+O(ar/n))+O((r/n)^2)\\
=&\,(1-\sigma e^{O(r^{4\eps}/n)})e^{\sigma}+O\left[\frac{r}{n}\sum_{a\ge 1}a\left(\frac{1}{(a-1)!}-\frac{1}{a!}\right)\right]+O((r/n)^2)\\
=&\,(1-\sigma e^{O(r^{4\eps}/n)})\,e^{\sigma e^{O(r^{4\eps}/n)}}+O(r/n)
=\,O(r/n).
\endalign
$$
This completes the proof of Theorem 2.1.\qed
\si

{\bf Note 2.1.\/}  Using Theorem 2.1 we will prove in the next section 3, Theorem 3.1, that it is quite
unlikely that the random digraph
with non-zero in/out-degrees has no isolated cycles and no sink/source-set with fewer
than $m/2$ arcs. 
In Section 4, we will use the  proof of Theorem 2.1 as a rough template for proving a genuinely sharp
asymptotic estimate for the probability of non-existence of ``simple sink/source-sets'' in 
the random digraph with {\it non-zero\/} in/out-degrees.   
This estimate coupled with Theorem 3.1 will deliver an asymptotic fraction of strongly-connected
digraphs among all such digraphs. 
\bi
{\bf 3. Bounding the number of the digraphs without simple sink-sets and small
complex sink-sets.\/}
\bi
 Let $\G11$ denote a digraph on $[n]$ which is chosen uniformly at
random from all such digraphs with $m$ arcs and with the smallest 
in-degree and the smallest out-degree both at least $1$. We call $S\subset [n]$
a sink-set (source-set resp.) if there is no arc $i\to j$ ($j\to i$ resp.) for
$i\in S$, $j\notin S$. A digraph is strongly-connected iff it has no
proper sink-sets and no proper source-sets. We call a sink-set (source-set resp.) $S$ simple,
if all the out-degrees (in-degrees) of $G_S$, the subdigraph induced by $S$, are
equal $1$, and complex otherwise. Our first step is to  the following 
result.
\si

\proclaim{Theorem 3.1} Suppose $r:=m-n\to\infty$ and $r=O(n)$.
Let $\Cal A_{n,m}$ denote the event ``$\G11$ has
a complex sink-set containing at most $m/2$ arcs, and has no simple
sink-set".  Then $\pr(\Cal A_{n,m})=O(n^{-1})$.
\endproclaim

\bi
{\bf Proof of Theorem  3.1.\/}  We will focus on  a core case $r=o(n)$, and at the end of the 
proof we will briefly discuss how to handle $r=\Theta(n)$.
\si 

Given $\nu\in [3,n)$, $\nu<\mu \le m_0$, let $X_{\nu,\mu}$
denote the total number partitions $[n]= A\uplus B$, ($|A|=\nu$),
such that (1) $A$ is a {\it minimal\/} complex sink-set with  
$\mu$ induced arcs, and (2) neither $G_A$ nor $G_B$ contain a simple
sink-set or source set. (Minimal means that $A$ does not contain a smaller
complex sink-set.) Let $\varepsilon>0$ be fixed. Set 
$$
\nu_0= \varepsilon\left(\frac{mn^2}{r^2}\right)^{1/3}\sim\varepsilon\,
\frac{n}{2r^{2/3}};\tag 3.1
$$
clearly $\nu_0=o(n)$. Define
$$
X=\sum_{\nu<\mu\atop \nu\le\nu_0,\mu\le m/2}X_{\nu,\mu}\,\,
+\sum_{\nu<\mu\atop\nu>\nu_0,\mu\le m/2}X_{\nu,\mu}.\tag 3.2
$$
In both sums $\nu$ and $\mu$ are subject to an additional condition, $\mu-\nu\le r$:
as we shall see shortly, $X_{\nu,\mu}=0$ otherwise.
The theorem will be proven when we show that $\ex[X]=O(n^{-1})$. 
\si

For this we will use the following bound: 
$$
\aligned
\ex[X_{\nu,\mu}]\le_b&\,\frac{r}{n}\,\,  E_{\nu,\mu},\\
 E_{\nu,\mu}:=&(x+1)\frac{\binom{n}{\nu}}
{\binom{m}{\mu}}\cdot\frac{\nu^{\mu-\nu}}{\mu\,(\mu-\nu)!}
\cdot \frac{(e^x-1)^{\nu}}{x^{\mu}}\cdot\frac{\la^{2\mu}e^{\nu\la}}
{(e^{\la}-1)^{2\nu}},
\endaligned\tag 3.3
$$
 uniformly for all $\nu<\mu$, $m-\mu\ge cm$, $c\in (0,1)$ being
fixed, and $x>0$. 
\si

To prove (3.3), let us bound the total number of partitions $(A,B)$ in question
in all the digraphs with parameters $n$, $m$, of minimum in-degree and
minimum out-degree $1$ at least. A set $A$ of cardinality $\nu$ can be chosen 
in $\binom{n}{\nu}$ ways. By
symmetry, we may consider $A=[\nu]$, and we will use $[n-\nu]$ to denote
$B=\{\nu+1,\dots,n\}$. Let us bound the total number of $G_A$'s.
Introduce $R$, the set consisting of $1$ and all vertices $j\in [n]$
reachable 
from $i$ by directed paths. $R\subseteq A$ as $A$ is a sink-set. Since $R$ is a sink-set
itself, $R=A$ by minimality of $A$. Thus there exists a directed tree $T$
rooted at $1$ that spans $A$, with arcs oriented away from the vertex $1$.
This spanning tree is such that the out-degree of $1$ equals its the out-degree
of $1$ in $G_A$ itself. Given 
$$
d_1\ge 1,\,\,d_2,\dots,d_{\nu}\ge 0,\quad (d_1-1)+d_2+\cdots+d_{\nu}=\nu-2,\tag 3.4
$$
there are 
$$
\frac{(\nu-2)!}{(d_1-1)!\prod_{i\ge 2}d_i!}
$$ 
rooted trees $T$ with outdegrees $d_1,\dots,d_{\nu}$.
The digraph $G_A$ is a disjoint union of $T$ and a complementary digraph $H$ on $[\nu]$. 
$\delta_i,\,\Delta_i$, ($i\in[\nu]$), respectively the in-degrees and the 
out-degrees of $H$, must meet the constraints
$$
\align
&\sum_i\delta_i=\sum_i\Delta_i=\mu-\nu+1,\tag 3.5\\
&d_1\ge 1,\,\delta_1\ge 1,\,\Delta_1=0;\quad \Delta_i+d_i\ge 1,\,(i\ge 2).\tag 3.6
\endalign
$$
(All in-degrees in $H$, except the root $1$, are positive, whence the single 
condition $\delta_1\ge 1$ ensures that $T\cup H$ has all its in-degrees positive.
The condition $\Delta_i+d_i\ge 1$, redundant for the root $1$, ensures that $T\cup H$
has all its out-degrees positive as well.) The total number of $H$'s with the in-degrees
$\delta_i$ and the out-degrees $\Delta_i$ is $(\mu-\nu+1)!/\prod_i(\delta_i!\Delta_i!)$,
at most. Hence the total number of $G_A$'s, with $A=[\nu]$, is on the order of
$$
\,(\nu-2)!(\mu-\nu+1)!\sum_{\bold d,\boldsymbol\delta,\boldsymbol\Delta}
\prod_i\frac{1}{(d_i-1_{\{i=1\}})!\,\delta_i!\,\Delta_i!},
$$
with $\bold d,\boldsymbol\delta,\boldsymbol\Delta$ meeting the constraints (3.4)-(3.6).
Using the generating functions,
$$
\align
&\sum_{\bold d,\boldsymbol\delta,\boldsymbol\Delta}
\prod_i\frac{1}{(d_i-1_{\{i=1\}})!\delta_i!\Delta_i!}\\
=&\,[x^{\nu-2}y^{\mu-\nu+1}z^{\mu-\nu+1}]
\sum_{\bold d,\boldsymbol\delta,\boldsymbol\Delta\,\atop\text{meet }(3.6)} 
\prod_i\frac{x^{d_i-1_{\{i=1\}}}y^{\delta_i}z^{\Delta_i}}{(d_i-1_{\{i=1\}})!\delta_i!\Delta_i!}\\
=&\,[y^{\mu-\nu+1}]\,(e^y-1)e^{(\nu-1)y}\\
&\times [x^{\nu-2}z^{\mu-\nu+1}]\left(\sum_{d\ge 1,\,\Delta= 0}\frac{x^{d-1}z^{\Delta}}
{d!\Delta!}\right)\cdot \left(\sum_{d+\Delta\ge 1}\frac{x^dz^{\Delta}}
{d!\Delta!}\right)^{\nu-1}\\
=&\,[y^{\mu-\nu+1}]\,(e^y-1)e^{(\nu-1)y}
\times [x^{\nu-2}z^{\mu-\nu+1}]\,e^x(e^{x+z}-1)^{\nu-1}.
\endalign
$$
Here
$$
[y^{\mu-\nu+1}]\,(e^y-1)e^{(\nu-1)y}\le [y^{\mu-\nu+1}]\, ye^{\nu y}=[y^{\mu-\nu}]e^{\nu y}
=\frac{\nu^{\mu-\nu}}{(\mu-\nu)!},
$$
and
$$
\align
[x^{\nu-2}z^{\mu-\nu+1}]\,e^x(e^{x+z}-1)^{\nu-1}\le&\,[x^{\nu-2}z^{\mu-\nu+1}]\,
e^{x+z}(e^{x+z}-1)^{\nu-1}\\
=&\,\frac{(\mu-1)!}{(\nu-2)!(\mu-\nu+1)!}\,[x^{\mu-1}]\,e^x(e^x-1)^{\nu-1},
\endalign
$$
where
$$
\align
[x^{\mu-1}]\,e^x(e^x-1)^{\nu-1}\le&\, \frac{e^x(e^x-1)^{\nu-1}}{x^{\mu-1}},
\quad (\forall\, x>0),\\
\le&\,(x+1)\frac{(e^x-1)^{\nu}}{x^{\mu}}.
\endalign
$$
Collecting the pieces, we bound the total number of $G_A$'s with $|A|=\nu$ by
$$
\multline
\binom{n}{\nu}\,(\nu-2)!(\mu-\nu+1)!\sum_{\bold d,\boldsymbol\delta,
\boldsymbol\Delta}
\prod_i\frac{1}{(d_i-1_{\{i=1\}})!\,\delta_i!\,\Delta_i!}\\
\le\,(x+1)\binom{n}{\nu}\,\frac{\nu^{\mu-\nu}\,(\mu-1)!}{(\mu-\nu)!}\,
\frac{(e^x-1)^{\nu}}{x^{\mu}}.
\endmultline\tag 3.7
$$
\si

Turn to the complementary subgraphs $G_B$, $B=[\nu]^c:=[n-\nu]$.
Let $\delta_i,\,\Delta_i$, $i\in [n-\nu]$, denote in-degrees and out-degrees of 
$G_{[n-\nu]}$, and let $\part_i$ denote the total number of arcs emanating
from $i\in [n-\nu]$ and ending at a vertex in $[\nu]$. Obviously
$$
\align
&\delta_i\ge 1;\quad \part_i+\Delta_i\ge 1,\quad i\in [n-\nu],\tag 3.8\\
&\sum_{i\in [n-\nu]}\delta_i=\sum_{i\in [n-\nu]}\Delta_i=m-\mu-\part.\tag 3.9
\endalign
$$
 Let us bound
the total number of the admissible $G_{[n-\nu]}$ for given $\part_i,\,i\in
[n-\nu]$. The second condition  in (3.8) simplifies to
$$
\Delta_i\ge 1,\quad i\in I(\boldsymbol\sigma):=\{i\in [n-\nu]\,:\,\part_i=0\}.
$$
Clearly, 
$$
|I(\boldsymbol\sigma)|\ge n-\nu -\part,\quad \part:=\sum_{i\in [n-\nu]}\part_i.
$$
By the first condition in (3.8), $n-\nu\le m-\mu-\part$; in particular, $\mu-\nu\le
m-n=r$. Actually $n-\nu
<m-\mu-\part$, since otherwise all $\delta_i=1$ and $G_{[n-\nu]}$ would
contain an isolated cycle, a special case of a simple sink-set. So 
$\part\le r-(\mu-\nu)< r$, and 
$$
|\,[n-\nu]\setminus I(\boldsymbol\sigma)|\le r\ll n-\nu,
$$
i.e. out-degrees $\Delta_i$ cannot be zero for all but at most $r$ 
specified vertices. 
Let $C_{\boldsymbol\sigma}(m-\mu-\part,n-\nu)$ denote the total number of digraphs
on $[n-\nu]$ with $m-\mu-\part$ arcs such that
$$
\delta_i\ge 1,\,\,(i\in [n-\nu]);\quad \Delta_i\ge 1,\,\,\,(i\in I(\boldsymbol\sigma)).  
\tag 3.10
$$
Further,
$$
C_{\boldsymbol\sigma}(m-\mu-\part,n-\nu)\le (m-\mu-\part)!\sum_{\boldsymbol{\delta}
,\boldsymbol{\Delta}
\text{ meet }(3.8), (3.9)}\,\prod_{i\in [n-\nu]}\frac{1}{\delta_i!\Delta_i!}.
$$
And, by the theorem 2.1, the total number of $G_{[n-\nu]}$ is on the order of
$$
\frac{m-\mu-(n-\nu -\part)}{n-\nu}\,C_{\boldsymbol\sigma}(m-\mu-\part,n-\nu)
\le_b\frac{r}{n}\,C_{\boldsymbol\sigma}(m-\mu-\part,n-\nu).
$$
Also, the total number of ways to choose in $[\nu]$
the partners of the vertices in $[n-\nu]$ is
$$
\prod_{i\in [n-\nu]}\binom{\nu}{\part_i}\le\nu^{\part}\!\!\prod_{i\in [n-\nu]}
\frac{1}{\part_i!}.
$$    
Therefore, for a given $\part$, the total number of the subdigraphs that 
complement $G_{[\nu]}$ is 
on the order of 
$$
\align
&\frac{r}{n}\,\nu^{\part}(m-\mu-\part)!\sum_{
\boldsymbol{\delta},\boldsymbol{\Delta}\text{ meet }(3.8), (3.9)}\,\,
\prod_{i\in [n-\nu]}\frac{1}{\part_i!\Delta_i!\delta_i!}\\
=&\,\frac{r}{n}\,\nu^{\part}(m-\mu-\part)!\,[u^{\part}v^{m-\mu-\part}]
\left(\sum_{i+j\ge 1}
\frac{u^i}{i!}\frac{v^j}{j!}\right)^{n-\nu}[w^{m-\mu-\part}]
\left(\sum_{k\ge 1}\frac{w^k}{k!}\right)^{n-\nu}\\
=&\,\frac{r}{n}\,\nu^{\part}(m-\mu-\part)!\,[u^{\part}v^{m-\mu-\part}]
(e^{u+v}-1)^{n-\nu}
\,[w^{m-\mu-\part}](e^w-1)^{n-\nu}\\
=&\,\frac{r}{n}\,\nu^{\part}(m-\mu-\part)!\,\binom{m-\mu}{\part}[u^{m-\mu}]
(e^u-1)^{n-\nu}
\,\,[w^{m-\mu-\part}](e^w-1)^{n-\nu}\\
\le_b&\,\frac{r}{n}\,(m-\mu)!\,\frac{\nu^{\part}}{\part!}\,((n-\nu)u)^{-1/2}
\frac{(e^u-1)^{n-\nu}}{u^{m-\mu}}
\cdot ((n-\nu)w)^{-1/2}\frac{(e^w-1)^{n-\nu}}{w^{m-\mu-\part}},
\endalign
$$
for all $u>0$, $w>0$. 

Setting $u=w=\la$, and summing for $\part\ge 0$,
we obtain an overall bound for the count of complementary subdigraphs:
$$
\multline
\frac{r}{n}\,(n\la)^{-1}(m-\mu)!\left(\sum_{\part\ge 0}\frac{(\nu\la)^{\part}}
{\part!}\right)\,\frac{(e^{\la}-1)^{2(n-\nu)}}
{\la^{2(m-\mu)}}\\
=\,\frac{r}{n}\,(n\la)^{-1}(m-\mu)!\,\frac{(e^{\la}-1)^{2(n-\nu)}e^{\nu\la}}{\la^{
2(m-\mu)}}.
\endmultline\tag 3.11
$$
(We have used $n-\nu=\Theta(n)$.)
\si

The total number of the partitions $(A,B)$ is bounded above by
the product of (3.7), (3.11) and $\binom{n}{\nu}$. Dividing this product by the
total number of the digraphs in question, we arrive at the bound (3.3). 
\bi

To get the most out of (3.3), we will use  $x=x(\nu,\mu)$, the 
minimum point of 
$$
h(\nu,\mu,x)=2\nu\ln(e^x-1)-2\mu\ln x,
$$
i.e. the root of
$$
h_x(\nu,\mu,x)=2\nu\frac{e^x}{e^x-1}-\frac{2\mu}{x}=0.\tag 3.12
$$
Considering $\nu$ and $\mu$ as continuously varying, we compute the
partial derivatives
$$
\aligned
&x_{\nu}=-\frac{\mu}{\nu^2}\,g(x),\quad x_{\mu}=\frac{1}{\nu}\,g(x),\\
&g(x)=\frac{(e^x-1)^2}{e^x(e^x-1-x)},
\endaligned\tag 3.13
$$
Since $g(0+)=2$ and $g(+\infty)=1$, it follows that $x(\nu,\mu)=\Theta((\mu-\nu)/\nu)$,
uniformly for all $0<\nu<\mu$. In particular, in (3.3) the factor $(x+1)/\mu$ is $O(1/\nu)$;
so we replace $x+1$ with $1/\nu$. 
(It is easy to show also that $x(\nu,\mu)>2(\mu-\nu)/\mu$.)
\bi

Consider  $\nu\le\nu_0$ and $\nu<\mu\le m/2$. By (3.3) and  (3.13), 
$$
\align
\frac{E_{\nu,\mu}}{E_{\nu,\mu-1}}\le &\,\la^2\frac{\mu}{m-\mu}\frac{\nu}{\mu-\nu}
\exp\bigl[h(\nu,\mu,x(\nu,\mu))-h(\nu,\mu-1,x(\nu,\mu-1))\bigr]\\
&(\tilde\mu\in [\mu-1,\mu],\quad \tilde x=x(\tilde\mu,\nu))\\
=&\,\la^2\frac{\mu}{m-\mu}\frac{\nu}{\mu-\nu}\exp\bigl[h_{\mu}(\nu,\tilde\mu,\tilde x)+
h_x(\nu,\tilde\mu,\tilde x)x_{\mu}(\nu,
\tilde\mu)\bigr]\\
=&\,\la^2\frac{\mu}{m-\mu}\frac{\nu}{\mu-\nu}\exp(-\ln \tilde x)\le 
\la^2\frac{\mu}{m-\mu}\frac{\nu}{(\mu-1)-\nu}\exp(-\ln \tilde x)\\
\le_b&\,\la^2\frac{\mu}{m-\mu}\,x^{-2}(\nu,\mu-1);\tag 3.14
\endalign
$$
($x(\nu,\mu)$ increases with $\mu$, $\frac{\mu-\nu}{\nu}=\phi(x(\nu,\mu))$ where
$$
\phi(x):=\frac{xe^x}{e^x-1}-1=\Theta(x),
$$
uniformly for $x>0$.) Let us show that the last expression decreases with $\mu$. 
Using the formula for $x_{\mu}$ in (3.13) for $\mu-1$ instead of $\mu$, 
we compute
$$
\align
&\frac{d}{d\mu}\left(\frac{\mu}{m-\mu}\,x^{-2}(\nu,\mu-1)\right)\\ 
=&\,\frac{m}{(m-\mu)^2x^2}-\frac{2}{(m-\mu)}\,\frac{\mu}{\nu}
\cdot\frac{(e^x-1)^2}{x^3e^x(e^x-1-x)}\qquad(x:=x(\nu,\mu-1))\\
\le&\,
\frac{m}{(m-\mu)^2x^2}-\frac{2}{(m-\mu)}\frac{xe^x}{e^x-1}
\cdot\frac{(e^x-1)^2}{x^3e^x(e^x-1-x)}\\
=&\,\frac{1}{x^2(m-\mu)^2}\left[m-2(m-\mu)\frac{e^x-1}{e^x-1-x}\right]<0,\quad (!)
\endalign
$$
as $\mu\le m/2$. Thus indeed the RHS in (3.14) decreases with $\mu$. At $\mu=\nu+2$,
for $\nu\le n_0=\varepsilon (mn^2/r^2)^{1/3}$, this RHS is
$$
\align
&\la^2\frac{\nu+2}{m-(\nu+2)}\,x^{-2}(\nu,\nu+1)\sim\,\la^2\frac{\nu}{m}
\left(\frac{2}{\nu}\right)^{-2}(1+O(1/\nu))\\
=&\,\frac{\la^2}{4m}\nu^3(1+O(1/\nu))
=\varepsilon^3\,\frac{4r^2/n^2}{4m}\frac{mn^2}{r^2}(1+O(1/\nu))<2\varepsilon^3,
\endalign
$$
if $\nu$ exceeds a large enough $\nu^*$. The same bound holds trivially for $\nu\le
\nu^*$, as $x(\nu,\nu+1)\ge 2/(\nu+1)$. Hence
$$
\frac{E_{\nu,\mu}}{E_{\nu,\mu-1}}\le_b 2\varepsilon^3\le \frac{1}{2},\quad (\nu\le n_0,\,\,\nu+2\le
\mu\le m/2),
$$
provided that $\varepsilon>0$ is chosen sufficiently small.
Consequently, denoting $x=x(\nu,\nu+1)$ and using
$$
\frac{\binom{n}{\nu}}{\binom{m}{\nu+1}}\le_b \frac{\nu}{m}\left(\frac{n}{m}\right)^{\nu},
\quad \frac{n}{m}=\frac{\la e^{\la}}{e^{\la}-1},\quad
\frac{e^x-1}{x}=e^x\frac{\nu}{\nu+1}=e^{2/\nu}(1+O(1/\nu))
$$
we bound
$$
\align
\frac{r}{n}\!\!\!\sum_{\nu<\mu\atop \nu\le n_0,\,\mu\le m/2}\!\!\!\!\!\!\!
\nu^{-1}E_{\nu,\mu}
\le&\,\frac{r}{n}
\sum_{\nu\le n_0}\nu^{-1}E_{\nu,\nu+1}\left(\sum_{j\ge 1}j\,(2/3)^{j-1}\right)\\
\le_b&\,\frac{r}{n}\sum_{\nu\le n_0}\frac{\binom{n}{\nu}}{\binom{m}{\nu+1}}
\cdot\frac{\la^2}{\nu x}\left[\frac{(e^x-1)^2}{x^2}\cdot\frac{\la^2e^{\la}}
{(e^{\la}-1)^2}\right]^{\nu}\\
\le_b&\,\frac{r}{n}\frac{\la^2}{m}\sum_{\nu\ge 1}\nu^2\left(\frac{\la}
{e^{\la}-1}\right)^{\nu}
\le_b \frac{r}{n}\frac{\la^2}{m}\left(1-\frac{\la}{e^{\la}-1}\right)^{-3}\\
=&\,O((r/n)(m\la)^{-1})=O(n^{-1}).
\endalign
$$
Thus
$$
\Sigma_1:=\frac{r}{n}\!\!\!\sum_{\nu<\mu\atop \nu\le n_0,\,\mu\le m/2}
\!\!\!\!\!\!\!\frac{1}
{\nu}\,E_{\nu,\mu}=O(n^{-1}). \tag 3.15
$$
\bi

It remains to bound  
$$
\Sigma_2:=\frac{r}{n}\!\!\!\sum_{\nu<\mu\atop \nu>n_0,\,\mu\le m/2}
\!\!\!\!\!\!\!\hat E_{\nu,\mu}.
$$
Since $k!=\Theta(k^{1/2}(k/e)^k)$, for $E_{\nu,\mu}$ in (3.3) we have
$$
E_{\nu,\mu}\le_b \frac{x+1}{\mu}\exp(H(\nu,\mu,x))\le_b\frac{1}{\nu_0}
\exp(H(\nu,\mu,x)),\tag 3.16
$$
where
$$
\aligned
H_{n,m}(\nu,\mu,x)=H(\nu,\mu,x):=&\,n\ln n-\nu\ln\nu-(n-\nu)\ln (n-\nu)\\
&-m\ln m+\mu\ln\mu+(m-\mu)\ln(m-\mu)\\
&+\nu\bigl(\ln(e^x-1)-2\ln(e^{\la}-1)+\la\bigr)\\
&+(\mu-\nu)\ln\nu-
(\mu-\nu)\ln\frac{\mu-\nu}{e}\\
&+\mu (2\ln\la-\ln x).
\endaligned\tag 3.17
$$
We will use $H_{\nu}$, $H_{\mu}$, $H_x$ to denote the partial derivatives of
$H(\nu,\mu,x)$. Like before, given $\nu<\mu$, we choose $x=x(\nu,\mu)$, the root of (3.12), or
equivalently the root of $H_x(\nu,\mu,x)=0$. Then
$$
\aligned
\frac{\partial H(\nu,\mu,x(\nu,\mu))}{\partial\nu}=&\,H_{\nu}(\nu,\mu,x(\nu,\mu))\\
=&\,-\ln\nu+\ln(n-\nu)+\ln\frac{e^{x(\nu,\mu)}-1}{(e^{\la}-1)^2}+\la\\
&+\frac{\mu-\nu}{\nu}-\ln\nu+\ln(\mu-\nu).
\endaligned\tag 3.18
$$
By (3.13), $x(\nu,\mu)$ strictly decreases with $\nu$, and then so does 
$\partial H(\nu,\mu,x(\nu,\mu))/\partial\nu$. That is, as a function of $\nu$,
$H(\nu,\mu,x(\nu,\mu))$ is {\it convex\/}. In fact, 
using $x_{\nu}(\nu,\mu)<0$ and $\mu-\nu\le r$, we have
$$
\aligned
&\frac{\partial^2 H(\nu,\mu,x(\nu,\mu))}{\partial\nu^2}=\,
\frac{\part H_{\nu}(\nu,\mu,x(\nu,\mu))}{\part\nu}\\
=&\,-\frac{1}{\nu}-\frac{1}{n-\nu}+\frac{e^{x(\nu,\mu)}x_{\nu}(\nu,\mu)}
{e^{x(\nu,\mu)}-1}-\frac{\mu}{\nu^2}-\frac{1}{\nu}-\frac{1}{\mu-\nu}\\
\le&-\frac{1}{r}.
\endaligned\tag 3.19
$$
\si

Since $x(\nu,\mu)\downarrow 0$
as $\nu\uparrow\mu$, and $x(\nu,\mu)\uparrow\infty$ as $\nu\downarrow 0$,
there exists a unique $\nu(\mu)$, the root of $\partial H(\nu,\mu,x(\nu,\mu))/
\partial\nu=0$, at which $H(\nu,\mu,x(\nu,\mu))$ attains its absolute maximum
$f(\mu):=H(\nu(\mu),\mu,x(\nu(\mu),\mu))$. Observe immediately that the
equation $H_{\nu}(\nu,\mu,x(\nu,\mu))=0$, with $H_{\nu}(\nu,\mu,x(\nu,\mu))$ given by 
(3.18), and (3.17) allow us to obtain a much simpler expression for $f(\mu)$,
namely
$$
\aligned
f(\mu)=&\,-n\ln(1-\nu/n)+\mu\ln\mu+(m-\mu)\ln(m-\mu)-m\ln m\\
&+\mu\ln\frac{\la^2\nu}{x(\mu-\nu)}.
\endaligned\tag 3.20
$$
Let us have a close look at $\nu(\mu)$ and
$x(\nu(\mu),\mu)$. 
\si

Observe first that, by (3.13),
$$
\frac{d x(\mu,\nu(\mu))}{d\mu}=x_{\mu}+x_{\nu}\nu_{\mu}=\frac{g(x)}{\nu}\left(1-
\frac{\mu}{\nu}\,\nu_{\mu}\right).\tag 3.21
$$
Differentiating $\partial H(\nu,\mu,x(\nu,\mu))/\partial\nu=0$ with respect
to $\mu$, using (3.13) and solving for $\nu_{\mu}$, we obtain
$$
\nu_{\mu}=\frac{(e^x-1)^{-1}g(x)/\nu}{2\nu^{-1}+(n-\nu)^{-1}+(\mu-\nu)^{-1}+\mu\nu^{-2}
+(e^x-1)^{-1}g(x)\mu\nu^{-2}}.
$$
So $\nu_{\mu}(\mu)\in (0,\nu(\mu)/\mu)$, i.e. $\nu(\mu)$ and, by (30), $x(\nu(\mu),\mu)$
both strictly increase with $\mu$.
\si
 
Turning the tables around, introduce the functions $\nu(x)$ and $\mu(x)$  
determined by two equations, (3.12) and $H_{\nu}(\nu,\mu,x)=0$, (see (3.18) for 
$H_{\nu}(\nu,\mu,x)$): 
$$
\nu=\nu(x):=n\,\frac{\dsize\frac{e^x-1}{(e^{\la}-1)^2}\phi(x)e^{\phi(x)}}
{1+\dsize\frac{e^x-1}{(e^{\la}-1)^2}\phi(x)e^{\phi(x)}},
\quad \mu=\mu(x):=\nu(x)\,
\frac{xe^x}{e^x-1},\tag 3.22
$$
where $\phi(x)=xe^x/(e^x-1)-1$.
Since $\nu(\mu)$, $x(\nu(\mu),\mu)$ are strictly increasing, then so are
$\nu(x)$ and $\mu(x)$. 
\si

Using the parameterization $\nu(x),\,\mu(x)$, let us evaluate
sharply $F(x):=f(\mu(x))\mathbreak= H(\nu(x),\mu(x),x)$ for $x\le c\lambda$, with
$c>0$ to be specified shortly.  After some algebra, (3.22) becomes
$$
\aligned
&\nu(x)=\,n\,\frac{0.5(x/\la)^2}{1+0.5(x/\la)^2}\,(1+\theta_1+O(\la^2)),\\
&\mu(x)=\,m\,\frac{0.5(x/\la)^2}{1+0.5(x/\la)^2}\,(1+\theta_2+O(\la^2)),\\
&\theta_1=\,\frac{7x/6-\la}{1+0.5(x/\la)^2},\quad
\theta_2=\,\theta_1+\frac{x-\la}{2};
\endaligned\tag 3.23
$$
(we have used $m/n=1+\la/2+O(\la^2)$). In particular, $\mu(x)$ is of order
$m\,(x/\la)^2$ exactly. Plugging these expressions into (3.20)
we obtain, after massive simplifications,
$$
\aligned
F(x)=&\,-r\,G(x/\la)+O(m x^2),\\
G(z):=&\,\ln\bigl(1+0.5z^2\bigr)-\frac{z^3}{1+0.5z^2}+
\frac{z^2}{1+0.5z^2},
\endaligned\tag 3.24
$$
($r=m-n$). $G(z)>0$ for $z\in (0,z_0]$, $z_0= 1.772$, and $G(z)\sim 1.5\,z^2$ as 
$z\downarrow 0$. For $x_0:=\la z_0$, it follows from (3.23) that
$$
m_1:=\mu(x_0)=m\,\frac{z_0^2}{1+z_0^2}(1+O(\la)) > 0.61 m.
$$
Returning to the pameterization $\nu(\mu),\,x(\nu(\mu),\mu)$, we see that 
$x(\nu(\mu),\mu)\le x_0$ for $\mu\le m_0=0.61m$,  and
$$
\align
f(\mu)=&\,F(x(\nu(\mu),\mu))\le -c_2r\left(\frac{x(\nu(\mu),\mu)}{\la}\right)^2+
O(mx^2)\\
=&\,-c_3\,r\,\frac{\mu}{m}+O(\mu\la^2)\le -c_4\,r\,\frac{\mu}{m},\tag 3.25
\endalign
$$
as $\la=O(r/m)$.
\si

With these preparations out of the way, turn directly to $\Sigma_2$. 
Given $\mu\in (\nu_0+1,m_1]$, ($n_0=\Theta(n/r^{2/3})$), by (3.16), (3.19) 
and (3.25),
$$
\align
\sum_{\nu\in (n_0,\mu)\atop\nu\ge \mu-r}E_{\nu,\mu}
\le_b&\,\nu_0^{-1}\exp\left(-c_4r\frac{\mu}{m}\right)
\times
\sum_{\nu\in (n_0,\mu)}
\exp\left(-c_1\,\frac{(\nu-\nu(\mu))^2}{r}\right)\\
\le_b&\,\nu_0^{-1}\exp\left(-c_4r\frac{\mu}{m}\right)\int_{-\infty}
^{\infty}e^{-c_1z^2/r}\,dz\\
\le_b&\,n^{-1}r^{7/6}\exp\left(-c_4r\frac{\mu}{m}\right).
\endalign
$$
Therefore
$$
\aligned
\frac{r}{n}\sum_{\mu\in (\nu_0,m_1]}
\sum_{\nu\in (n_0,\mu)\atop\nu\ge \mu-r} E_{\nu,\mu}
\le_b&\,\frac{r^{13/6}}{n^2}\sum_{\mu>\nu_0}\exp\left(-c_4r\frac{\mu}{m}\right)\\
\le_b&\,\frac{r^{13/6}}{n^2}\frac{\exp(-c_4r\nu_0/m)}{1-e^{-c_4r/m}}\\
\le_b&\,\frac{r^{7/6}}{n}\,e^{-c_4r^{1/3}}=O(1/n).
\endaligned\tag 3.26
$$
\si
Combining (3.26) and (3.15), we conclude that $\ex[X]=O(n^{-1})$. This completes
the proof of Theorem 3.1 for the case $r=o(n)$. Since  $r/n\to 0$ however slowly, we actually 
established our claim for $r\to\infty$ such that $r\le\eps n$, where $\eps>0$ is sufficiently small.
\bi

Turn to $r=\Theta(n)$.   Up to and including the equations (3.22) the proof is 
basically the same. Furthermore, for $x\le x_0$, $x_0$ being small enough, $F_{n,m}(x)=H(\nu(x),\mu(x),
x)=
-n\Theta(x^2)$. The challenge is to show that $F_{n,m}(x)=-\Theta(n)$ for all $x\ge x_0$, such that
$\mu(x)\le m/2$. This is definitely so for $r\le\eps n$. Let  
$$
c^*=\sup\bigr\{c\ge \eps\,:\, F_{n,n+cn}(x)< 0\text{ if } x>x_0 \text{ and } \mu(x)\le m/2\bigl\}.
$$
If $c^*<\infty$, then for $m=n+c^*n$ there exists $x^*>x_0$ such that $\mu(x^*)\le m/2$ and
$F_{n,m}(x^*)=0$, $F_{n,m}^\prime(x^*)=0$. Equivalently, we must have $f_{n,m}(\mu)=0$, $f^\prime_{n,m}(\mu)=0$ for some $0<\mu\le m/2$, with $f_{n,m}(\mu)=f(\mu)$ defined by (3.20).
\si

Since $H_x(\mu,\nu,x(\nu,\mu))=0$, $H_{\nu}(\nu(\mu),\mu,x(\nu,\mu))=0$, by (3.17) we have
$$
f^\prime_{n,m}(\mu)=H_{\mu}(\mu,\nu(\mu),x(\nu(\mu),\mu))=\ln\left(\frac{\mu}{m-\mu}\,\frac{\nu}{\mu-\nu}\,
\frac{\la^2}{x}\right)=0,\tag 3.27
$$
so that, denoting $x=x(\nu(\mu),\mu)$,
$$
\mu=\frac{m}{1+\dsize\frac{\la^2}{x\phi(x)}},\quad\phi(x)=\frac{x e^x}{e^x-1}-1.\tag 3.28
$$
Since $\mu\le m/2$, it must be true that $\la^2/x\phi(x)\ge 1$, so that $\la >\phi(x)$, as $\phi(x)<x$ for $x>0$.
Furthermore, combining (3.27) and (3.20), we obtain
$$
f_{n,m}(\mu)=-n\ln(1-\nu(\mu)/n)+m\ln(1-\mu/m)=0.\quad (!)\tag 3.29
$$
Combining (3.28)-(3.29) with (3.22), and denoting $\psi(x)=x e^x/(e^x-1)$, we see that 
there must exist a solution $(\la,x)$ of the following system of two equations:
$$
\aligned
&1+\frac{e^x-1}{(e^{\la}-1)^2}\,\phi(x)e^{\phi(x)}=\left(1+\frac{x\phi(x)}{\la^2}\right)^{\psi(\lambda)},\\
&\frac{\psi(\la)}{\psi(x)}\,\frac{1}{1+\dsize\frac{\la^2}{x\phi(x)}}=
\frac{\dsize\frac{e^x-1}{(e^{\la}-1)^2}\phi(x)e^{\phi(x)}}
{1+\dsize\frac{e^x-1}{(e^{\la}-1)^2}\phi(x)e^{\phi(x)}}.
\endaligned\tag 3.30
$$
Our  task is to show that assuming existence of such a solution $(\la,x)$ we get
 a contradiction.  

Since $\psi(\la)>1$, the first equation in (3.30) implies that
$$
1+\frac{e^x-1}{(e^{\la}-1)^2}\,\phi(x)e^{\phi(x)}>1+\frac{x\phi(x)}{\la^2}\psi(\la),\tag 3.31
$$
or, using the definition of $\psi(\cdot)$,
$$
\frac{e^x-1}{x}\,e^{\phi(x)}>\frac{e^{\la}-1}{\la}\, e^{\la}.
$$
As $\phi(x)<\la$, and $(e^y-1)/y$ is increasing, it follows that $x>\la$. By (3.31), the second equation
yields an inequality
$$
\frac{\psi(\la)}{\psi(x)}\,\frac{1}{1+\dsize\frac{\la^2}{x\phi(x)}}>\frac{\dsize\frac{x\phi(x)}{\la^2}\psi(\la)}
{1+\dsize\frac{x\phi(x)}{\la^2}\psi(\la)},
$$
or
$$
\psi(x)\left(1+\frac{\la^2}{x\phi(x)}\right)<\frac{\la^2}{x\phi(x)}+\psi(\la),
$$
or
$$
\frac{\la^2}{x\phi(x)}\bigl(\psi(x)-1\bigr)<\psi(\la)-\psi(x).
$$
This is impossible since $\psi(x)>1$ and $\psi(\la)<\psi(x)$.
\si

Thus, for $r=\Theta(n)$, $F_{n,m}(x)=-n\Theta(x^2)$ for $x\le x_0$, and $F_{n,m}(x)=-\Theta(n)$,
for all $x\ge x_0$ such that $\mu(x)\le m/2$. 

With this property established, the bound for $\Sigma_2$ can be proved in essentially
the same was as for $r=o(n)$. \qed
\bi

{\bf Note 3.2.\/} If a digraph with positive in/out degrees has a sink-set (a source-set)
$T$ with $|T|<n$, then there is a source-set (a sink-set) $\Cal T\subseteq T^c$. (If, say,
$T$ is a source-set, then, for a vertex $v\in T^c$, the set consisting of $v$ and all vertices
reachable from $v$ is a sink-set in $T^c$.) Therefore if a digraph of this sort is not
strongly-connected {\it and\/} has no simple sink/source-sets, it must have either a complex
sink-set or a complex source-set with at most $m/2$ induced arcs. By Theorem 3.1, the
probability of this event is $O(n^{-1})$. Thus, it remains to find a sharp estimate
for the probability of no simple sink/source-sets and to check that this probability far
exceeds $n^{-1}$.

\bi
{\bf 4. Counting the digraphs without simple sink/source-sets.\/}  
\bi
A subgraph induced by a simple sink-set (source-set resp.) is a disjoint union
of cycles.  So our task in this section is to find a sharp asymptotic formula for
the probability that $G_{1,1}(n,m)$ has no induced cycles $\Cal C$ such that
there are no arcs $i\to j$ ($j\to i$ resp.) with $i\in \Cal C$ and $j\notin \Cal C$. 
\si

Consider first the corresponding probability for the random multi-digraph
\linebreak $MG_{1,1}(n,m)$ induced by the $m$-long insertion sequence chosen
uniformly at random among all $h_{1,1}(n,m)$ such sequences. (See Note 1.2
and (1.17) for the definition of $h_{1,1}(n,m)$.) Let $X=X_n$ denote the total
length of the cycles in all simple sink-sets and source-sets.
\proclaim{Lemma 4.1}
$$
\aligned
\ex\bigl[z^X\bigr]=&\,\frac{(m!)^2}{h_{1,1}(n,m)}\sum_a\frac{(n)_a}{(m)_a}\\
&\times [(x_1x_2)^{m-a}y^a]\biggl(H(\bold x,y,z)\prod_{i=1}^2(e^{x_i}-1)^{n-a}\,
\biggr),
\endaligned\tag 4.1
$$
where
$$
H(\bold x,y,z)=\frac{1-yz}{1-y}\prod_{i=1}^2\frac{1-ye^{x_i}}
{1-yze^{x_i}}.\tag 4.2
$$
In particular,
$$
\multline
\pr(X=0)=\ex\bigl[0^X\bigr]=\,\frac{(m!)^2}{h_{1,1}(n,m)}\sum_a\frac{(n)_a}{(m)_a}\\
\times [(x_1x_2)^{m-a}y^a]\left((1-y)^{-1}
\prod_{i=1}^2(1-ye^{x_i})
(e^{x_i}-1)^{n-a}\right).
\endmultline\tag 4.3
$$
\endproclaim
{\bf Proof of Lemma 4.1.\/} Let us first derive a formula for $\ex\bigl[\binom{X}{k}\bigr]$.
$\binom{X}{k}$ is the total number of unordered $k$ tuples of vertices belonging to the
cycles of simple sink/source-sets.  By symmetry, $\ex\bigl[\binom{X}{k}\bigr]$ equals 
$\binom{n}{k}$ times the probability that the cycles contain vertices from $[k]$. To evaluate
this probability we need to count the number of $m$-long insertion sequences that result in 
multi-digraphs in which the set $[k]$ belongs to those
cycles. We find it a bit easier to turn things around. First, we evaluate the number of the insertion
sequences  having
several such (disjoint) cycles induced by a subset of vertices. Second, we multiply this number by the total number
of ways to select a set of $k$ points contained in these cycles, under a constraint that {\it every\/} cycle
is represented in such a set.  Third, we sum these products over all those collections of disjoint
sets. Finally we divide the sum by $h_{1,1}(n,m)$, the total number of the insertion sequences.
\si

 Let $a$ denote a (generic) cardinality of a chosen set $A$ of vertices; there are $\binom{n}{a}$
 such sets. Let  $a_{1,1},\dots,
 a_{1,\ell_1}$
($a_{2,1},\dots, a_{2,\ell_2}$, resp.) denote the lengths of sink (source, resp.) cycles 
to be formed out of the
$a$ vertices. As loops allowed,  $a_{i,j}\ge 1$ and
$$
\sum_{i=1}^2\sum_{j=1}^{\ell_i}a_{i,j}=a.\tag 4.4
$$
The total number of ways to form $\ell_1+\ell_2$ directed cycles is
$$
\frac{a!}{\ell_1!\ell_2!}\prod_{i,j}\frac{(a_{i,j}-1)!}{a_{i,j}!}=
\frac{a!}{\ell_1!\ell_2!}\prod_{i,j}\frac{1}{a_{i,j}}.
\tag 4.5
$$
The total number of ways to select $k$ points from these cycles subject to the representation 
constraint is
$$
\align
\sum_{k_{i,j}\ge 1\,:\,\|\bold k\|=k}\prod_{i,j}\binom{a_{i,j}}{k_{i,j}}=&\,[w^k]\prod_{i,j}
\left(\sum_{\kappa\ge 1}\binom{a_{i,j}}{\kappa}w^{\kappa}\right)\\
=&\,[w^k]\prod_{i,j}\bigl[(1+w)^{a_{i,j}}-1\bigr].\tag 4.6
\endalign
$$
The product of the counts in (4.5) and (4.6) is
$$
\frac{a!}{\ell_1!\ell_2!}\prod_{i,j}\frac{(1+w)^{a_{i,j}}-1}{a_{i,j}}.\tag 4.7
$$
Given these cycles $C_{1,1},\dots, C_{1,\ell_1},C_{2,1},\dots,
C_{2,\ell_2}$, let $\hat\delta_i$, $1\le i\le \ell_1$, denote the number of arcs $u\to v$ 
with $u\notin C_{1,i}$ and $v\in C_{1,i}$, and let $\hat\Delta_i$, $1\le i\le \ell_2$,
denote the number of arcs $u\to v$ with $u\in C_{2,i}$ and $v\notin C_{2,i}$. Since
$C_{1,i}$ is a sink cycle, there are no arcs from $C_{1,i}$ to its complement;
likewise, there are no arcs leading to $C_{2,i}$ from its complement. To avoid
overcounting we consider an isolated cycle only as a source cycle. Thus we have
only the constraint $\hat\delta_i\ge 1$, $1\le i\le \ell_1$. For $i\notin A$, let
$\delta_i\ge 1$, $\Delta_i\ge 1$ denote its in-degree and its out-degree. Clearly
$$
\sum_{i=1}^{\ell_1}\hat\delta_i+\sum_{i\notin A}\delta_i=m-a,
\quad \sum_{i=1}^{\ell_2}\hat\Delta_i+\sum_{i\notin A}\Delta_i=m-a.\tag 4.8
$$
By (1.5), the total number of the $(m-a)$-long insertion sequences that result
in a multi-digraph on the vertex set $\{1,\dots,\ell_1,\ell_1+1,\dots,\ell_1+\ell_2, [n]
\setminus A\}$, with the in-degrees $\hat\delta_1,\dots,\hat\delta_{\ell_1}$,
$\delta_i$, ($i\in [n]\setminus A$), and the out-degrees $\hat\Delta_1,\dots,
\hat\Delta_{\ell_2}$, $\Delta_i$, ($i\in [n]\setminus A$), is
$$
((m-a)!)^2\prod_{i=1}^{\ell_1}\frac{1}{\hat\delta_i!}\,
\prod_{i=1}^{\ell_2}\frac{1}{\hat\Delta_i!}\,\prod_{i\in [n]\setminus A}\frac{1}
{\delta_i!\Delta_i!}.
$$
We need to multiply this count by
$$
\prod_{i=1}^{\ell_1}a_{1,i}^{\hat\delta_i}\cdot\prod_{i=1}^{\ell_2}a_{2,i}^{\hat\Delta_i},
$$
since, say, for each of $\hat\delta_i$ arcs ending at the ``aggregated'' vertex, that
represents the cycle $C_{1,i}$, there are $a_{1,i}=|C_{1,i}|$ ways to specify the arc
end. There is another factor still missing, namely $(m)_{m-a}$, which counts
the total number of ways to select those $(m-a)$ positions in the whole $m$-long
insertion sequence when a newly added arc belongs to that multi-digraph, rather than to one
of the cycles.  By (4.8), and the constraints on $\hat\delta_i$, $\delta_i,\Delta_i$,
the resulting product equals
$$
\aligned
m!(m-a)!\,&[x_1^{m-a}]\,(e^{x_1}-1)^{n-a}\prod_{i=1}^{\ell_1}(e^{a_{1,i}x_1}-1)\\
\times\,&[x_2^{m-a}]\,(e^{x_2}-1)^{n-a}\prod_{i=1}^{\ell_2}e^{a_{2,i}x_2}\,.
\endaligned\tag 4.9
$$
Further, the product of the factors in (4.7) and (4.9) dependent on $a_{1,i},\,a_{2,i}$,
equals
$$
\aligned
&\frac{1}{\ell_1!}\,\prod_{i=1}^{\ell_1}\frac{(e^{a_{1,i}x_1}-1)\bigl[(1+w)^{a_{1,i}}-1\bigr]}{a_{1,i}}\\
\times\,&
\frac{1}{\ell_2!}\,
\prod_{i=1}^{\ell_2}\frac{e^{a_{2,i}x_1}\bigl[(1+w)^{a_{2,i}}-1\bigr]}{a_{1,i}}.
\endaligned\tag 4.10
$$
Summing  the products (4.10) over $a_{i,j}\ge 2$ meeting the constraint (4.4), and
then over all $\ell_1$, $\ell_2$ we get
$$
\aligned
&[y^a]\,\sum_{\ell_1\ge 0}\frac{1}{\ell_1!}\left(\sum_{b\ge 1}
y^{b}\,\frac{(e^{b x_1}-1)\bigl[
(1+w)^{b}-1\bigr]}{b}\right)^{\ell_1}\\
&\times\,
\sum_{\ell_2\ge 0}\frac{1}{\ell_2!}\left(\sum_{b\ge 1}y^{b}\,\frac{e^{b x_2}\bigl[
(1+w)^{b}-1\bigr]}{b}\right)^{\ell_2}\\
=&\,[y^{a}]\,\exp\left[\sum_{b\ge 1}y^{b}\,\frac{(e^{b x_1}-1)\bigl[
(1+w)^{b}-1\bigr]}{b}+
\sum_{b\ge 1}y^{b}\,\frac{e^{bx_2}\bigl[
(1+w)^{b}-1\bigr]}{b}\right]\\
=&\,[y^{a}]\,H(\bold x,y,1+w),
\endaligned\tag 4.11
$$
see (4.2) for the definition of  $H(\bold x,y,z)$. 
\si

Collecting all the pieces, we write
$$
\ex\left[\binom{X}{k}\right]=\frac{(m!)^2}{h(n,m)}\,[w^k]\,\sum_a\frac{(n)_a}{(m)_a}
\,[(x_1x_2)^{m-a}]\,H(\bold x,y,1+w)\prod_{i=1}^2(e^{x_i}-1)^{n-a}.
$$
Consequently
$$
\multline
\ex\bigl[(1+w)^X\bigr]=\sum_kw^k\,\ex\left[\binom{X}{k}\right]\\
=\frac{(m!)^2}{h(n,m)}\,\sum_a\frac{(n)_a}{(m)_a}
\,[(x_1x_2)^{m-a}]\left(H(\bold x,y,1+w)\prod_{i=1}^2(e^{x_i}-1)^{n-a}\right),
\endmultline
$$
which is equivalent to (4.1). \qed
\bi

Armed with this lemma, we obtain an asymptotic formula for $\pr(X=0)$.
\proclaim{Theorem 4.2} For the random multi-digraph $MG_{1,1}(n,m)$, with
$r=m-n\to\infty$, $m=O(n)$,
$$
\pr(X=0)=\frac{\left(1-\dsize\frac{\la}{e^{\la}-1}\right)^2}{
1-\dsize\frac{\la}{e^{\la}(e^{\la}-1)}}\,\bigl(1+O(r^{-1}\ln^2 r)\bigr),\tag 4.12
$$
with $\la$ determined by $\la e^{\la}/(e^{\la}-1)=m/n$.
\endproclaim
{\bf Proof of Theorem 4.2.\/} In the formula (4.3) for $\pr(X=0)$, 
$$
(1-y)^{-1}\prod_{i=1}^2(1-ye^{x_i})=1+\prod_{i=1}^2(e^{x_i}-1)\,\sum_{a\neq 0,1}y^a-\bigl[1+
\sum_{i=1}^2(e^{x_i}-1)\bigr]\,y
$$
and
$$
h_{1,1}(n,m)=(m!)^2\bigl([x^m](e^x-1)^n\bigr)^2.
$$
So (4.3) becomes
$$
\pr(X=0)=1+\sum_{1\le \a< n}\frac{(n)_a}{(m)_a}Q^2(a,a-1)-\frac{n}{m}\bigl(Q(1,1)+Q(1,0)\bigr)^2,
\tag 4.13
$$
where
$$
Q(a,b):=\frac{[x^{m-a}]\,(e^x-1)^{n-b}}{[x^{m}]\,(e^x-1)^{n}}.
$$
 As in the proof of the theorem 2.1, for $a,b<n$,
$$
Q(a,b)\le_b\frac{\la^a}{(e^{\la}-1)^b}\cdot\sqrt{\frac{n}{n-b}},\quad \frac{(n)_a}{(m)_a}\le
\left(\frac{n}{m}\right)^a=\left(\frac{e^{\la}-1}{\la e^{\la}}\right)^a,
$$
and consequently, setting
$$
A_n:=\left\lceil\frac{2\ln r}{\la}\right\rceil=\Theta\bigl(r^{-1}n\ln r\bigr),
$$
after simple estimates we have 
$$
\sum_{a\ge A_n}\frac{(n)_a}{(m)_a}\,Q^2(a,a-1)\le_b\,e^{-A_n\la}\min\{1,\la\}\le r^{-2}\min\{1,\la\}.\tag 4.14
$$
(We compare the series tail with $\min\{1,\la\}$ because the latter is the order of the RHS in (4.12).)
\si

Turn to $a< A_n$. Let us have a close look at $[x^{m-a}]\,(e^x-1)^{n+1-a}$,
the numerator in the fraction $Q(a,a-1)$. Introduce $\tilde\la$, satisfying
$$
\frac{x e^x}{e^x-1}=\frac{m}{n+1},
$$
and the independent positive Poissons $\tilde Y^i=Y^i(\tilde\la)$. Using (2.7)
we have
$$
\multline
[x^{m-a}]\,(e^x-1)^{n+1-a}=\,\frac{(e^{\tilde\la}-1)^{n+1-a}}{\tilde\la^{m-a}}
\,\,\pr\!\left(\sum_{i=1}^{n+1-a}\tilde Y^i=m-a\right)\\
=\,\frac{(e^{\tilde\la}-1)^{n+1-a}}{\tilde\la^{m-a}}\,\,
\pr\!\left(\sum_{i=1}^{n+1}\tilde Y^i=m\right)\,\bigl(1+O(an^{-1}+a^2rn^{-2})
\bigr)
\endmultline
$$
A bit of calculus, based on (3.13), shows that $\tilde\la=\la+O(n^{-1})$.
Consequently
$$
\frac{(e^{\tilde\la}-1)^{n+1-a}}{\tilde\la^{m-a}}=
\frac{(e^{\la}-1)^{n+1-a}}{\la^{m-a}}\,\bigl(1+O(r^{-1}+n^{-1}a)\bigr),
$$
and
$$
\pr\left(\sum_{i=1}^{n+1}\tilde Y^i=m\right)=\frac{1+O(r^{-1})}
{(2\pi (n+1)\var[\tilde Y])^{1/2}}=\bigl(1+O(r^{-1})\bigr)
\pr\left(\sum_{i=1}^{n}Y^i=m\right).
$$
Therefore
$$
\align
Q(a,a-1)=&\,\bigl(1+O(r^{-1}+an^{-1}+a^2rn^{-2})\bigr)
\frac{\la^a}{(e^{\la}-1)^{a-1}}\\
=&\bigl(1+O(r^{-1}\ln^2 r)\bigr)\frac{\la^a}{(e^{\la}-1)^{a-1}},\tag 4.15
\endalign
$$
uniformly for $a\le A_n$.
\si

Also, uniformly for  $a<A_n$,
$$
\frac{(n)_a}{(m)_a}=\left(\frac{n}{m}\right)^a\,\bigl(1+O(a^2rn^{-2})\bigr)=
\left(\frac{n}{m}\right)^a\,\bigl(1+O(r^{-1}\ln^2r)\bigr).\tag 4.16
$$
Using (4.15)-(4.16),  $m/n=\la e^{\la}/(e^{\la}-1)$,  and denoting $\sigma=\la/[e^{\la}(e^{\la}-1)]$,
we compute:
$$
\multline
\sum_{1\le a<A_n}\frac{(n)_a}{(m)_a}\,Q^2(a,a-1)=(e^{\la}-1)^2\left(\sum_{1\le a<A_n}\sigma^a\right)
\bigl(1+O(r^{-1}\ln^2 r)\bigr)\\
=\,\frac{\sigma (e^{\la}-1)^2}{1-\sigma}\bigl[1+O(\sigma^{A_n})+O(r^{-1}\ln^2 r)\bigr]\\
=\,\frac{\sigma (e^{\la}-1)^2}{1-\sigma}\bigl(1+O(r^{-1}\ln^2 r)\bigr),
\endmultline\tag 4.17
$$
as
$$
\sigma^{A_n}\le e^{-(1-\sigma)A_n}\le e^{-\la A_n}\le r^{-1}.
$$
For the last, negative,  term in (4.13) we need to be more precise. Picking a simple contour $L$ enclosing the origin, by Cauchy integral formula we have
$$
\align
[x^m]\,(e^x-1)^n=&\frac{1}{2\pi i}\oint\limits_L\frac{(e^z-1)^n}{z^{m+1}}\,dz\\
=&\,\frac{1}{2\pi i}\cdot\frac{n}{m}\oint\limits_L\frac{(e^z-1)^{n-1}e^z}{z^m}\,dz\\
=&\,\frac{n}{m}\bigl\{[x^{m-1}]\,(e^x-1)^n+[x^{m-1}]\,(e^x-1)^{n-1}\bigr\}.
\endalign
$$
So, recalling the definition of $Q(\cdot,\cdot)$,
$$
Q(1,1)+Q(1,0)=\frac{m}{n}.\quad (!)\tag 4.18
$$
Combining (4.14), (4.16) and (4.17), we easily transform (4.13)  into (4.12).
\qed
\bi

Turn now to the random digraph $G_{1,1}(n,m)$. Let $P_{n,m}$ denote the probability that 
$G_{1,1}(n,m)$ has no simple sink-sets and source sets.  The theorem 4.2 provides the
asymptotic estimate for $\Cal P_{n,m}=\pr(X=0)$, the corresponding probability for the random multigraph
$MG_{1,1}(n,m)$.  We can write
$$
\Cal P_{n,m}=\frac{|S^*|}{h_{1,1}(n,m)};
$$
here $S^*$ is the set  of all $m$-long sequences of arc-insertions that result in a multi-digraph
without the simple sink/source sets. Denoting by $\boldsymbol\delta(\bold s)$ and $\boldsymbol\Delta(\bold s)$
the in/out degrees of a multi-digraph corresponding to a generic sequence $\bold s$, we see that the total
number of the {\it graphic\/} sequences $\bold s$ leading to a digraph without simple
sink/source sets is
$$
h_{1,1}^G(n,m):=\sum_{\bold s\in S^*}F(\boldsymbol\delta(\bold s),\boldsymbol\Delta(\bold s)).
$$
By (1.22), and $\Cal P_{n,m}=\Theta(r/n)$,
$$
\align
h_{1,1}^G(n,m)=&\,\bigl(e^{-\eta}+O(r^{\eps}n^{-1/2+\gamma})\bigr)|S^*|+O\bigl(h_{1,1}(n,m)e^{-r^{c\eps}
n^{2\gamma}}\bigr)\\
=&\,e^{-\eta}\,\Cal P_{n,m}h_{1,1}(n,m)\bigl(1+O(r^{\eps}n^{-1/2+\gamma})+
O(ne^{-r^{c\eps}n^{2\gamma}})\bigr)\\
=&\,e^{-\eta}\Cal P_{n,m}h_{1,1}(n,m)\bigl(1+O(r^{\eps}n^{-1/2+\gamma})\bigr).
\endalign
$$
Therefore
$$
P_{n,m}=\frac{h^G_{1,1}(n,m)}{h_{1,1}(n,m)}=e^{-\eta}\Cal P_{n,m}\bigl(1+O(r^{\eps}n^{-1/2+\gamma})\bigr),
\tag 4.19
$$
for $\eps>0$ and $\gamma\in (0,1/2)$.
\si

Combining (4.19), (4.12) and recalling the formula (1.12) for $\eta$, we have proved the following.
\proclaim{Theorem 4.3} For $r\to\infty$ and $m=O(n)$,
$$
\align
P_{n,m}=&\,\frac{\left(1-\dsize\frac{\la}{e^{\la}-1}\right)^2}{
1-\dsize\frac{\la}{e^{\la}(e^{\la}-1)}}\,\exp\left(-\frac{m}{n}-\frac{\la^2}{2}\right)\\
&\times \bigl(1+O(r^{-1}\ln^2 r+r^{\eps}n^{-1/2+\gamma})\bigr),\tag 4.20
\endalign
$$
for $\eps>0$ and $\gamma\in (0,1/2)$.
\endproclaim
\bi

By (4.20), $P_{n,m}=\Theta(\la)=\Theta(r/n)\gg n^{-1}$. In light of Note 3.2, (4.20)
combined with Theorem 1.2 complete the proof of our claim (1) in the
introduction.
\bi
{\bf Acknowledgement.\/} Back in $2001$, during my memorable visit to Nick Wormald at
University of Melbourne, we spent most of our time together working on the alternative
derivation of the  Bender-Canfield-McKay formula for the count of connected
sparse graphs, the remarkable formula that struck us by how sharp it was. It occurred to us then 
that the embedding idea might also
work for counting the strongly-connected directed graphs. It was immediately
clear though that an intrinsically harder notion of strong connectivity presented
a difficult new challenge. In Summer of $2009$  Nick informed me that
together with Xavier P\'erez they had made a further progress toward obtaining
a fully-proved asymptotic formula for the counts of strongly-connected digraphs .
In particular, this formula would yield a directed counterpart of the BCM formula, but without
an explicit estimate of an error term.  Reenergized  by the news, I set up to
see if I could  obtain a sharp asymptotic formula, with
a remainder term qualitatively matching that of the BCM formula. The approach 
in this paper is  naturally very different from the one used by P\'erez and Wormald.
However the paper would not be possible without
the joint work with Nick on the BCM formula, and our initial attacks on the directed
case. I owe to Nick my debt of gratitude.

\enddocument